\documentclass[brochure,12pt,french]{bourbaki}

\usepackage[T1]{fontenc}
\usepackage{lmodern,amssymb,bm}
\usepackage{graphicx}

\usepackage[french]{babel}

\usepackage[utf8]{inputenc}
\usepackage{microtype}
\usepackage{float}

\usepackage[
backend=bibtex,
style=authoryear, 
citestyle=authoryear-comp
]{biblatex}

\AtEveryBibitem{%
  \clearfield{issn} 
  \clearfield{doi} 

  \ifentrytype{online}{}{
    \clearfield{url}
  }
}
\DeclareFieldFormat[misc]{title}{\<#1\>} 

\usepackage{csquotes}
\newcommand{\N}{\mathbb {N}}
\newcommand{\R}{\mathbb {R}}

\newcommand{\T}{\mathbb {T}}
\newcommand{\D}{\mathbb {D}}
\newcommand{\Z}{\mathbb {Z}}

\newcommand{\Ic}{\mathcal {I}}
\newcommand{\Jc}{\mathcal {J}}
\newcommand{\Mc}{\mathcal {M}}
\newcommand{\Pc}{\mathcal {P}}
\newcommand{\Uc}{\mathcal {U}}
\newcommand{\Ss}{\mathbb {S}}
\newcommand{\Gc}{\mathcal {G}}
\newcommand{\cF}{\mathcal {F}}
\newcommand{\cI}{\mathcal {I}}

\theoremstyle{plain} 
\newtheorem{nota}[defi]{Notation}
\newtheorem*{remo}{Remords}
\newtheorem*{ques}{Question}

\newtheorem*{quesprovis}{Questions provisoires}

\date{Novembre 2019}
\bbkannee{72\textsuperscript{e} année, 2019--2020}  
\bbknumero{1166}                                      

\title{La démonstration de la  conjecture de l'entropie positive d'Herman }
\subtitle{d'après Berger et Turaev}

\author{Marie-Claude Arnaud}
\address{IMJ-PRG, UMR7586\\
Université de Paris \\
Bâtiment Sophie Germain,\\
 Boite Courrier 7012\\
8 Place Aurélie Nemours,\\
75205 Paris Cedex 13
}
\email{marie-claude.arnaud@imj-prg.fr}

\begin{document}

\maketitle

\noindent{\bf Abstract}  We present the proof of Berger and Turaev of Herman's positive entropy conjecture. In every neighbourhood of identity in the set of smooth symplectic diffeomorphisms of the 2-dimensional disc, there exists a diffeomorphism with positive metric entropy.

\section{Introduction}

Ce texte traite de l'entropie métrique pour la mesure de Lebesgue en dimension 2, une quantité qui mesure le chaos détecté par la mesure de Lebesgue pour un système dynamique différentiable la préservant.  Bien  longtemps, les seuls exemples connus d'entropie métrique positive\footnote{Ici, positive signifie strictement positive.}  ont été les systèmes Anosov (appelés aussi uniformément hyperboliques), définis sur le tore 2-dimensionnel. D'autres exemples ont suivi dans les années 80, construits par  \textcite{Katok1980},   sur toute surface, mais toujours élaborés à partir des Anosov, et toujours loin de l'identité au sens de la topologie~$C^1$.
Au congrès international de 1998,   Herman énonce la conjecture suivante.
\begin{conj}[Herman]\label{conjherman}
Soit ${\rm Diff}^\infty_\omega(\D)$ l'ensemble des difféomorphismes du disque de classe $C^\infty$ qui préservent l'aire. Dans tout voisinage de l'identité dans ${\rm Diff}^\infty_\omega(\D)$, il existe  un difféomorphisme d'entropie métrique positive.
\end{conj}

Une autre conjecture fameuse de l'entropie positive, qui concerne la famille standard des difféomorphismes du tore définis par $f_\lambda (x, y)=(2x - y + \lambda\sin 2\pi x, x)$, a été énoncée par \textcite{Sinai1994}.
\begin{conj}[Sinai]\label{conjsinai}
Pour tout paramètre $\lambda\not=0$, l'entropie métrique de $f_\lambda$ est positive.
 \end{conj}
 En ce qui concerne la conjecture~\ref{conjsinai}, on ne sait démontrer pour aucun paramètre $\lambda$ que l'entropie est positive, même si les simulations numériques semblent montrer des mers  de chaos. 
 
  La conjecture~\ref{conjherman} vient d'être démontrée par \textcite{BergerTurave2019}, dans un  véritable tour de force et l'objet de ce texte et de présenter leur démonstration.

 Après avoir rappelé ce qu'est l'entropie métrique, nous commencerons par décrire ce que Berger et Turaev appellent les îlots stochastiques, qui généralisent un exemple dû à  \textcite{Przytycki1982}. Ces objets, porteurs d'entropie positive, sont fragiles et ne résistent pas aux perturbations. Berger et Turaev construisent   un îlot stochastique qui résiste à un certain type de perturbations qu'ils introduisent, dites relatives aux liens.  Nous parlerons alors de renormalisation et expliquerons comment l'existence de bandes homoclines pour un point périodique hyperbolique permet de créer par perturbation des dynamiques qui, renormalisées, sont très proches de n'importe quelle dynamique (par exemple de l'îlot stochastique de Berger et Turaev). Une partie des méthodes utilisées ici ont été développées par  \textcite{Turaev2003}, dans le contexte des dynamiques universelles et par 
 \textcite{GonchenkoShilnikovTuraev2008}, pour étudier les domaines de Newhouse. Il s'agit de la partie la plus technique de la preuve. Finalement, nous montrerons comment construire dans tout voisinage de l'identité  un difféomorphisme $f$ ayant une bande homocline, ce qui permettra d'injecter une dynamique proche de l'îlot stochastique de Berger et Turaev comme renormalisée de $f$. Il faudra ensuite restaurer l'îlot stochastique. Cette dernière étape utilise un opérateur astucieusement défini pour évaluer l'écart entre deux branches de variétés invariantes.

 \begin{remo}
   L'élément central de l'exposé est la démonstration de la conjecture de
   Herman. Mais en vérité, Berger et Turaev démontrent plus que cette
   conjecture et montrent le résultat suivant~: soit $f$  un difféomorphisme préservant
   l'aire qui a un point périodique non hyperbolique et $V$ un voisinage de $f$ en topologie  $C^\infty$. Alors $V$ contient un  difféomorphisme préservant
   l'aire d'entropie métrique positive. \\
   Ceci leur permet de progresser vers la
   conjecture suivante.
   \begin{conj}[Berger \& Turaev]\label{conjbergerturaev}
     Il existe un ensemble dense de ${\rm Diff}^\infty_\omega(\D)$ dont tout
     élément est d'entropie métrique positive.
   \end{conj}
 \end{remo}
 
 Je remercie chaleureusement Sylvain Maillot pour sa relecture de cet exposé et ses précieux conseils, ainsi que Pierre Berger et Sylvain Crovisier qui m'ont  aidé à améliorer tant le texte que les illustrations. Merci aussi à toute l'équipe éditoriale qui fait un travail remarquable dans un délai très resserré.
\section{Parlons d'entropies}
Un {\em système dynamique discret} est la donnée d'une transformation (bijection) $f:X\rightarrow X$ d'un ensemble $X$ muni d'une structure dont on étudie les itérées $(f^n)_{n\in\N}$. En particulier, 
\begin{itemize}
\item si $X$ est un espace métrique, $f$ est un homéomorphisme;
\item si $X$ est un espace de probabilité, $f$ est mesurable et préserve la probabilité.
\end{itemize}
Pour introduire une notion de complexité d'un tel système, on s'intéresse aux
\og morceaux finis\fg{} d'orbites $(f^ix)_{0\leq i\leq n}\in X^{n+1}$ de longueur $n\in\N$ et on regarde ce qui se passe quand $n$ tend vers $+\infty$.  Étant donné un recouvrement  ouvert ou mesurable   
$$X=\bigcup_{0\leq i\leq m}X_i$$
de $X$, chaque morceau fini d'orbite $(f^ix)_{0\leq i\leq n}$ a au moins  un itinéraire $(n_i(x))_{0\leq i\leq n}\in [0, m]^{n+1}$ qui vérifie $f^i(x)\in X_{n_i(x)}$. On note $\Ic((X_i), f,  n)\subset [0, m]^{n+1}$ l'ensemble de tous les itinéraires pour $f$. Si les $X_i$ forment une partition et pas seulement un recouvrement, chaque morceau d'orbite a un unique itinéraire $I(x)$ et l'entropie de $f$   en temps $n$ relativement à la partition $(X_i)$ est alors le logarithme du nombre des itinéraires.
$$h((X_i), f, n)=\log |\Ic((X_i), f, n)|\footnote{Les barres $|.|$ servent à désigner le cardinal.}$$
Si les $X_i$ forment seulement un recouvrement, on choisit $\Jc((X_i), f, n)\subset \Ic((X_i),f,  n)$ de cardinal minimal tel que chaque morceau d'orbite de longueur $n$ a un itinéraire dans $\Jc((X_i),f, n)$, et on note
$$h((X_i), f, n)=\log |\Jc((X_i), f, n)|.$$
Dans les deux cas, la quantité $h((X_i), f, n)$ est sous-additive en $n$ et donc la quantité
$$h((X_i), f)=\lim_{n\rightarrow \infty}\frac{1}{n}h((X_i), f, n)$$
existe. C'est l'entropie de $f$ relativement au recouvrement  $(X_i)$, qui mesure donc le taux de croissance exponentiel du nombre d'itinéraires en fonction de $n$. On a ainsi
$$\Jc((X_i),f,  n)\varpropto \exp\left( nh((X_i), f)\right).\footnote{Le signe $\varpropto$ signifie que le logarithme du rapport des deux quantités est négligeable devant $n$.}$$

Quand $X$ est un espace métrique et $f$ un homéomorphisme, l'{\em entropie
  topologique} $h_{\rm top}(f)$ de $f$ est alors le supremum des $h((X_i), f)$
pour $(X_i)$ parcourant l'ensemble des recouvrements ouverts de $X$. Cette
définition d'entropie topologique est due à 
\textcite{AdlerKonheimMcAndrew1965} qui s'inspiraient de la définition d'entropie métrique de Kolmogorov que nous donnons plus loin. Une définition équivalente souvent plus maniable pour traiter les exemples fut introduite indépendamment par \textcite{Dinaburg1971} et 
\textcite{Bowen1970} par la suite.
 \begin{exem}
Une isométrie d'un espace métrique compact est d'entropie topologique nulle.
\end{exem}

Quand $(X,\mu)$ est un espace de probabilité et $f$ est une bijection mesurable qui préserve $\mu$, on définit l'entropie métrique en tirant parti de la mesure.  
Pour $(X_i)_{0\leq  i\leq m}$ partition mesurable, pour chaque $I\in \Ic((X_i), f, n)$, on note $X_I$ l'ensemble des points de~$X$ d'itinéraire $I$. L'entropie métrique  de $f$   en temps $n$ relativement à la partition~$(X_i)$ est alors 
$$h((X_i), f, \mu, n)=-\sum_{I\in\Ic}\mu(X_I)\log\left( \mu(X_I)\right).$$
C'est l'intégrale de la fonction $x\in X\mapsto -\log\mu(X_{I(x)})$. On a   alors
$$h((X_i), f, \mu, n)\leq\log |\Ic((X_i), f, n)|=h((X_i), f, n).\footnote{C'est une conséquence de la convexité de la fonction $x\mapsto x\cdot \log(x)$.}$$
avec égalité   seulement   dans le cas particulier où les $X_I$ sont équiprobables. L'entropie métrique de $f$ relativement à la partition $(X_i)$ est  
$$h((X_i), f, \mu)=\lim_{n\rightarrow \infty} \frac{1}{n}h((X_i), f, \mu, n).$$
L'{\em entropie métrique} $h_\mu(f)$ de $f$ est  le supremum des
$h((X_i),\mu,  f)$ pour $(X_i)$ parcourant l'ensemble des partitions
mesurables de $X$. Cette notion, antérieure à la notion d'entropie
topologique, fut introduite par Kolmogorov, en 1959
(cf. \cite{ Kolmogorov1958}), qui s'inspirait des travaux de Shannon concernant
la théorie de l'information.

On suppose que $X$ est un espace métrique compact et que $f$ en est un
homéomorphisme. L'ensemble $\Mc(X)$ des probabilités boréliennes de $X$ est
muni de la topologie faible et l'ensemble $\Mc(f)$ des probabilités
boréliennes invariantes par $f$ est alors un compact non vide. 
\textcite{Goodwyn1969}  a montré  l'inégalité variationnelle
$$\forall \mu\in \Mc(f), h_\mu(f)\leq h_{\rm top} (f).$$
Il est donc possible que l'entropie d'une mesure particulière soit nulle et que l'entropie topologique soit positive\footnote{Les deux notions d'entropie sont toujours positives ou nulles. Pour ne pas alourdir le texte et en suivant la terminologie anglaise, nous avons omis le terme strictement devant positive, mais dans ce texte positive signifiera toujours strictement positive.}\footnote{Pour un difféomorphisme de classe $C^1$ d'une variété compacte, les deux types d'entropie sont finies.}, alors que l'inverse est impossible.
 \textcite{Katok1980}  a   donné une condition nécessaire et suffisante  pour qu'un difféomorphisme de classe $C^{1+\alpha}$ d'une surface soit d'entropie topologique positive~: un tel difféomorphisme a un fer à cheval\footnote{Le fer à cheval est un ensemble compact invariant  qui apparait quand les variétés stable et instable d'un point périodique hyperbolique s'intersectent transversalement ailleurs qu'en ce point. Ils furent construits par  \textcite{Smale1965}.}. Or, il y a de très nombreux difféomorphismes qui ont un fer à  cheval. Dans un travail en cours, Le Calvez et Sambarino montrent que si $M$ est une surface compacte munie d'une forme d'aire, pour tout $k\in[1, \infty]$, il existe un ouvert dense de l'ensemble des difféomorphismes de classe $C^k$ de $M$ qui préservent l'aire  dont tout élément a un fer à  cheval.  Précisons qu'un fer à cheval est un ensemble de Cantor\footnote{Rappelons qu'un ensemble de Cantor est un ensemble non vide, compact, totalement discontinu et sans point isolé. Il a donc la puissance du continu.} et est donc petit au sens de la topologie car c'est un fermé d'intérieur vide. 

Plaçons nous donc dans le cas d'une surface compacte munie d'une forme d'aire, par exemple le disque unité de $\R^2$ muni de la mesure de Lebesgue. Désormais, nous ne nous intéresserons à l'entropie métrique que pour la mesure associée à la forme d'aire. On sait alors que l'entropie topologique de la plupart (au sens de la catégorie de Baire) des difféomorphismes du disque qui préservent l'aire est positive. Ceci n'implique bien sûr pas que leur entropie métrique (pour la mesure de Lebesgue donc) est positive.

Connaissons-nous beaucoup de difféomorphismes d'entropie métrique positive? En 1977,  Pesin a fait le lien entre l'entropie métrique pour une mesure de Lebesgue et le taux exponentiel avec lequel des points proches s'éloignent sous l'action dynamique, i.e. les taux  de divergence des différentielles $Df^n$.  Ceux-ci s'appellent les exposants de Lyapunov. Plus précisément, pour un difféomorphisme $f$ qui préserve une mesure de probabilité $\mu$ équivalente à une forme volume, le fibré tangent peut se décomposer comme somme directe $\oplus E_i(x)$ de fibrés mesurables $Df$-invariants sur un ensemble de $\mu$-mesure~$1$ et pour chaque $v\in E_i(x)$ non nul, la limite suivante existe et est indépendante de $v$. 
$$\lambda_i(x)=\lim_{n\rightarrow +\infty}\frac{1}{n}\log\left( \| Df^n(x)v\|\right) .$$
Il s'agit des  exposants de Lyapunov au point $x$.

\textcite{Pesin1977}  montre que si $f$ est  un difféomorphisme de classe au moins $C^2$ d'une variété qui préserve une mesure $\mu$ équivalente à une forme volume, alors  l'entropie est égale à la somme des exposants de Lyapunov positifs, i.e.{}
$$h_\mu(f)=\int \sum\lambda_i^+(x){\rm dim }E_i(x)d\mu(x)$$
où les $\lambda_i^+$ sont les exposants de Lyapunov positifs et  $\dim E_i(x)$ est la multiplicité de~$\lambda_i^+(x)$. Ceci permet de construire des exemples de difféomorphismes d'entropie métrique positive pour la mesure de Lebesgue.
\begin{exem} 
La transformation du tore 2-dimensionnel $f:\T^2=\left(\R/\Z\right)^2\rightarrow \T^2$ définie par $f(\theta_1, \theta_2)=(2\theta_1+\theta_2, \theta_1+\theta_2)$ a pour entropie  $
{\log\left( \frac{3+\sqrt{5}}{2}\right)}$ puisque les valeurs propres de la matrice $\left(\begin{smallmatrix} 2&1\\
 1&1\end{smallmatrix}\right)$ sont $
{\frac{3\pm\sqrt{5}}{2}}$. 
\end{exem}
Cet exemple appartient à la classe plus vaste des difféomorphismes Anosov, ceux pour lesquels on peut écrire le tangent à la variété comme somme directe de deux sous-fibrés $Df$-invariants $E^s\oplus E^u$, l'un suivant lequel $Df$ contracte et l'autre suivant lequel $Df$ dilate. On sait que la seule surface compacte qui porte un difféomorphisme Anosov est le tore $\T^2$ et qu'un difféomorphisme Anosov n'est pas homotope à l'identité, et donc est loin de l'identité. On sait aussi qu'être Anosov est une propriété stable par perturbation~$C^1$ petite, et on obtient donc un ouvert de difféomorphismes dont l'entropie métrique est positive. Nous n'avons aucune raison de penser que les exemples présentant de l'entropie métrique positive que nous allons construire par la suite ont un voisinage en topologie~$C^\infty$ qui a  la même propriété, mais nous ne savons démontrer ni la véracité de ceci ni son contraire.

Si on n'a  en tête que les difféomorphismes Anosov comme exemples ayant une entropie métrique positive, plusieurs  questions paraissent naturelles.

\begin{quesprovis}
  \footnote{Provisoires car comme nous allons l'expliquer dans ce texte, Katok
    a répondu positivement à la première et Berger et Turaev positivement à la seconde.}{\sl
    \begin{enumerate}
    \item Existe-t-il des exemples de difféomorphismes d'entropie métrique
      positive sur des surfaces autres que le tore?
    \item Peut-on trouver de tels difféomorphismes très proches de l'identité
      en topologie~$C^\infty$?
    \item Plus généralement, l'ensemble des difféomorphismes d'entropie
      métrique positive est-il dense dans l'ensemble des difféomorphismes qui
      préservent l'aire?
    \end{enumerate}}
\end{quesprovis}
\section{Les îlots stochastiques et leurs perturbations}\label{secile}
\subsection{Notion et construction d'îlots stochastiques}\label{secexilot}
\textcite{Katok1979}  a prouvé que toute surface compacte porte au moins un difféomorphisme de classe $C^\infty$ qui préserve la forme d'aire et est d'entropie métrique positive. Il commence par  construire un tel difféomorphisme pour le disque qui vaut l'identité au bord, puis injecte ce difféomorphisme dans toute surface. 

Pour construire son exemple sur le disque, il part d'un difféomorphisme $A$  Anosov du tore tel que $A\circ( -{\rm Id})=-A$ et tel que les quatre points fixes de $-{\rm Id}$ sont aussi fixés par $A$. Il ralentit $A$ aux points fixes (i.e. le rend tangent à l'identité) en préservant à la fois la symétrie par $-{\rm Id}$ et l'invariance des directions stables et instables. Il obtient ainsi un difféomorphisme $F$ du tore tel que $F\circ( -{\rm Id})=-F$ et qui en dehors des quatre points fixes a un exposant de Lyapunov positif. Le quotient du tore par $-{\rm Id}$ est une sphère. En passant donc au quotient, on obtient un homéomorphisme~$f$ de~$\Ss^2$. Le passage au quotient définit un revêtement de degré 2 ramifié aux quatre points fixes, $f$ est de classe $C^\infty$ en dehors des points fixes. Le seul problème de régularité qui pourrait se poser concernerait les points fixes, mais avoir rendu le difféomorphisme~$F$  tangent à l'identité permet d'obtenir un $f$ de classe $C^\infty$, qui a comme avait $F$ une entropie métrique positive. Ensuite, on éclate la sphère en un des points fixes pour obtenir finalement un difféomorphisme $g$ du disque fermé qui vaut l'identité au bord et qui est d'entropie métrique positive.

\begin{nota}
  On notera $\D$ le disque unité de $\R^2$ et ${\rm Diff}^r_\omega(\D)$
  désignera l'ensemble des difféomorphismes symplectiques de classe $C^r$ du
  disque unité.
\end{nota}

 Remarquons que dans cet exemple, presque tout point (pour la mesure d'aire) a
 un exposant de Lyapunov positif, i.e. le domaine qui porte l'entropie
 positive est de mesure pleine.  Cette propriété  ne persiste pas par
 perturbation petite en topologie~$C^\infty$~: on peut rendre générique\footnote{Un point fixe est elliptique si sa
   différentielle est conjuguée à une rotation. Il est alors générique si le
   rapport de l'angle de cette rotation à $2\pi$  est diophantien  et le
   premier terme dans son développement en forme normale est non nul.} l'un des points
 fixes elliptiques, de
 telle sorte qu'on puisse appliquer en son voisinage les théorèmes K.A.M. qui
 fournissent un ensemble de mesure positive sur lequel tout point a ses
 exposants de Lyapunov qui sont nuls\footnote{Cet ensemble s'écrit comme la
   réunion d'un ensemble de courbes sur lesquelles la dynamique est conjuguée
   à une rotation.  Cette propriété d'avoir   des courbes invariantes KAM
   subsiste d'ailleurs si on fait des perturbations petites en topologie
   $C^\infty$, même si le fait d'avoir un nombre de rotation diophantien ne
   subsiste pas.}. Mais il pourrait subsister une région invariante plus
 petite de mesure non nulle sur laquelle un des exposants de Lyapunov est
 positif, la question est ouverte. Tout comme la question suivante plus
 précise d'Herman à  l'ICM de 1998 \parencite{Herman1998}.

 \begin{ques} {\sl Comme l'application qui à $f$ associe son plus grand
     exposant de Lyapunov pour la mesure de Lebesgue est semi-continue
     supérieurement, l'ensemble $\Gc$ des difféomorphismes de classe $C^1$
     préservant l'aire dont les exposants de Lyapunov sont nuls presque
     partout est un ensemble $G_\delta$. Par la théorie de Baire, on ne peut
     avoir qu'une des deux situations suivantes.
 \begin{enumerate}
 \item soit $\Gc$ est dense;
 \item soit il existe un ensemble ouvert en topologie~$C^1$ dont tout élément est d'entropie métrique positive.
 \end{enumerate}
 Quelle est la situation qui est vérifiée, la première ou la seconde?}
\end{ques}

 En topologie~$C^1$, \textcite{Mane1983,Mane1996}
 avait annoncé 
puis donné des éléments de démonstration 
qu'un difféomorphisme symplectique générique d'une surface qui n'est pas Anosov est d'entropie métrique nulle, la démonstration complète ayant été donnée par  \textcite{Bochi2002}. La question reste ouverte pour la topologie~$C^\infty$.
 
 \textcite{Przytycki1982}  a construit sur le tore un exemple où on a une situation mixte~: une région où l'entropie métrique est positive et une région   elliptique où l'entropie métrique est nulle. Berger et Turaev appellent {\em îlots stochastiques} une famille de dynamiques que nous décrivons ci-après en définition~\ref{defilot} et  qui englobe la région d'entropie métrique positive construite par Przytycki.

Soit $f$ un difféomorphisme de classe $C^\infty$ d'une surface $M$ qui préserve une forme d'aire. On rappelle qu'un point $\tau$-périodique $P=f^\tau(P)$ est hyperbolique si les deux valeurs propres $\lambda$, $\lambda^{-1}$ de $Df^\tau(P)$ sont réelles et de module différent de $1$.  On sait alors qu'il existe deux immersions injectives $j_s, j_u:(\R, 0)\hookrightarrow (M, P)$ de classe $C^\infty$ telles que $j_s(\R)$ (resp. $j_u(\R)$) est la variété stable (resp. instable) de $P$ qui est définie par
$$W^s(P, f)=\{ x\in M; \lim_{n\rightarrow +\infty} d(f^nx, f^nP)=0\}=W^u(P, f^{-1}).$$
Une {\em branche stable} (resp. {\em instable}) de $P$ est la réunion de $\{ P\}$ et d'une des deux composantes connexes de $W^s(f, P)\backslash \{ P\} $ (resp. $W^u(f, P)\backslash \{ P\}$). 
\begin{defi}
On appelle alors {\em lien hétérocline} un cercle topologiquement plongé (i.e. un lacet) invariant qui s'écrit comme une union finie de branches stables de points périodiques hyperboliques. Si $n$ points périodiques interviennent, on parle de {\em $n$-lien}, de {\em bi-lien} si $n=2$.
\end{defi}
\begin{figure}[H]
\centering\includegraphics[width=4.5cm]{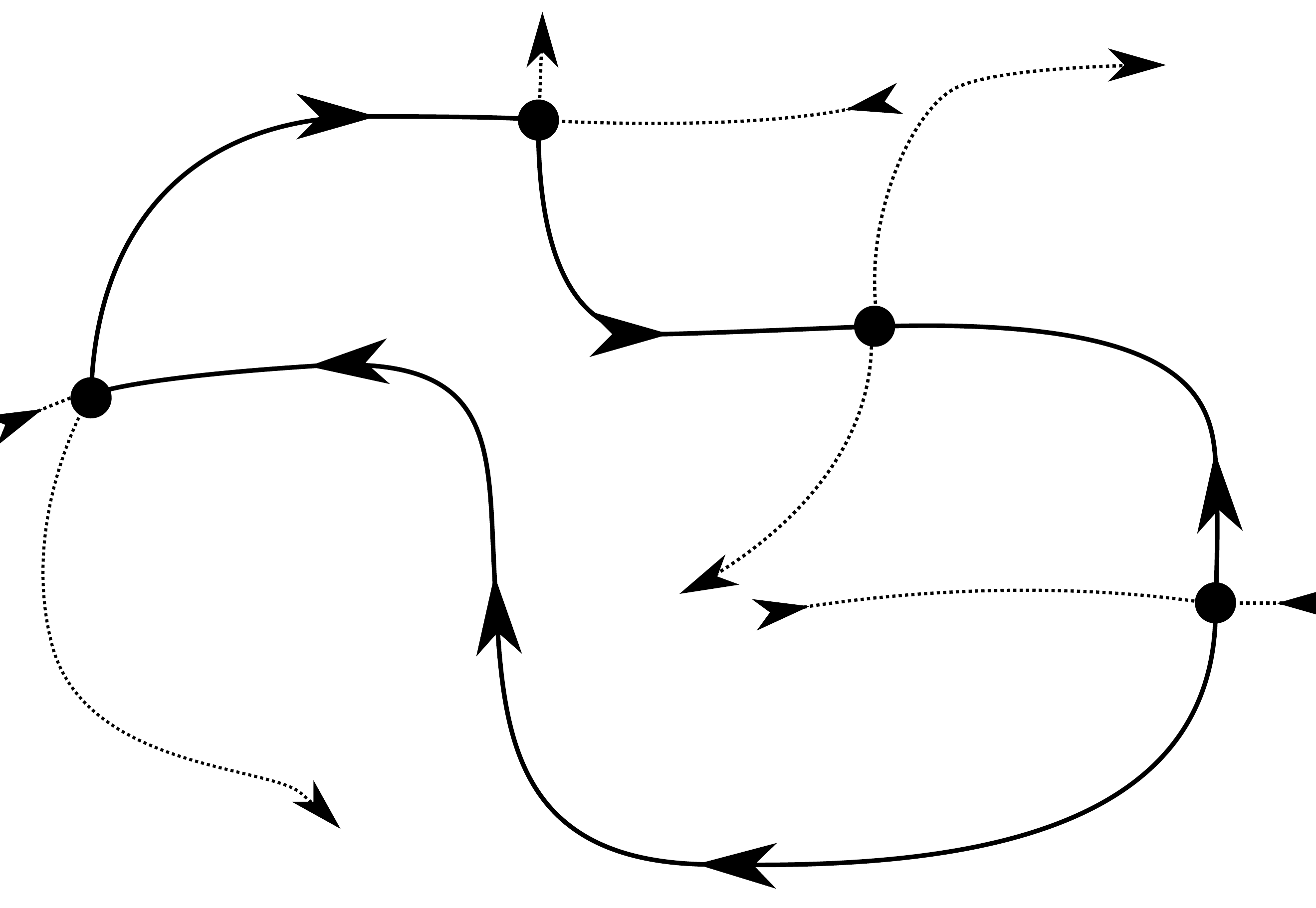}
\caption{Un lien hétérocline en gras}
\end{figure}

Un lien hétérocline est alors de classe $C^\infty$ sauf éventuellement aux points périodiques où il peut présenter un point anguleux.
\begin{defi}\label{defilot}
  Un  îlot stochastique est un ouvert invariant $\cI$ dont le bord est une union finie de liens hétéroclines tel que chaque point de $\cI$ a un exposant de Lyapunov positif.
\end{defi}

\textcite{BergerTurave2019} adaptent la construction de Katok en remplaçant les quatre points fixes par des disques bordés par des bi-liens hétéroclines, le complémentaire de leur union étant alors un îlot stochastique $\Ic_0$. Ils se ramènent au disque comme faisait Katok.

\begin{figure}[H]
\centering\includegraphics[width=4cm]{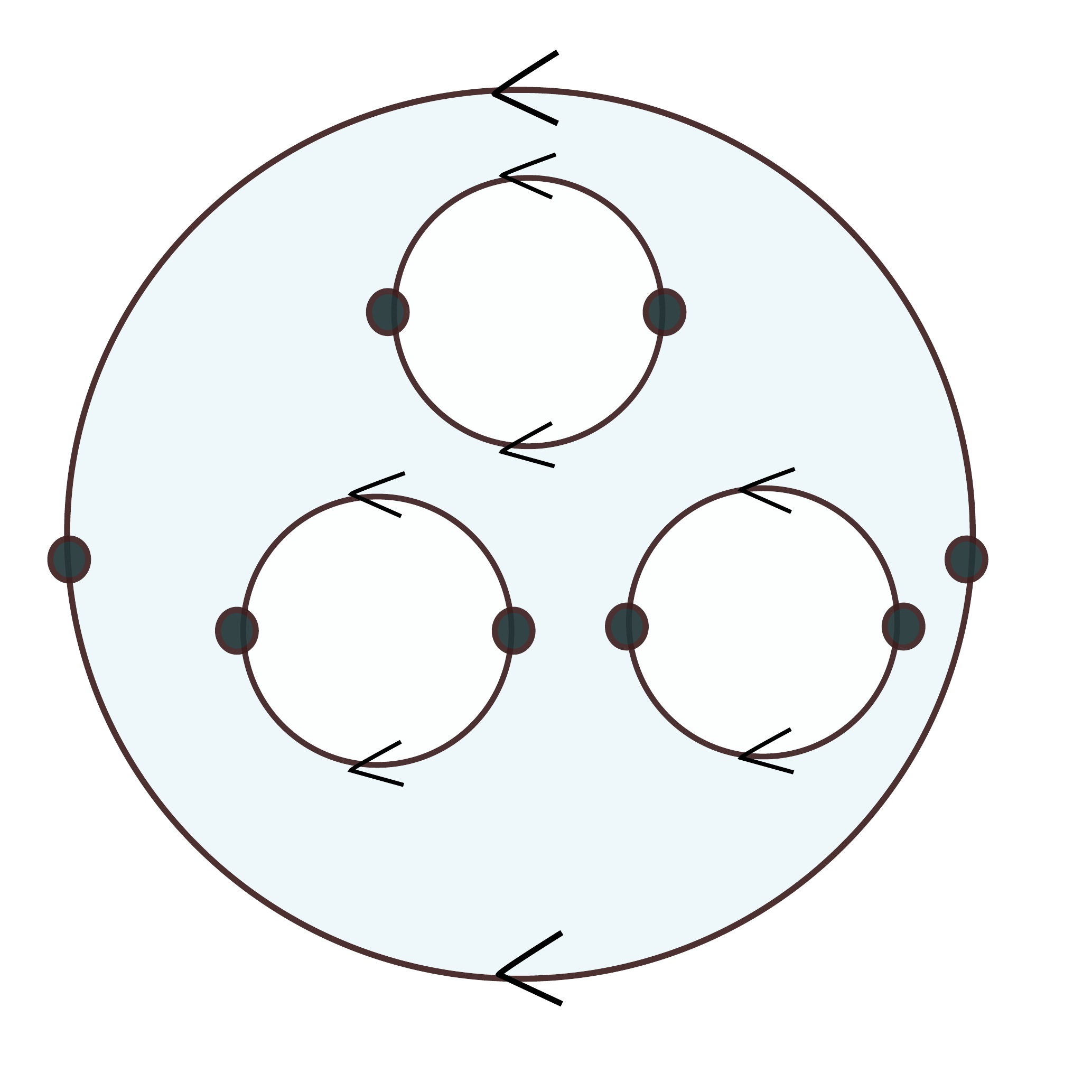}
\caption{L'îlot stochastique de Berger et Turaev.}
\end{figure}
\subsection{Ilots stochastiques et perturbations}
Comme nous l'avons déjà mentionné, le problème de ces constructions est leur fragilité. Par perturbation, les points hyperboliques et leurs branches persistent en ce sens qu'ils admettent des prolongement par continuité, mais  en général, deux branches ne coïncident pas, et quand elles se coupent, elle se coupent transversalement. Or, pour injecter leur construction dans une dynamique proche de l'identité (ce que n'est pas cet exemple), les auteurs vont avoir besoin de perturber. Ils proposent alors de restreindre l'ensemble des perturbations de manière à ne considérer que des perturbations qui préservent les liens hétéroclines (plutôt d'ailleurs ils feront des perturbations qui cassent ces liens mais seront capables de les restaurer).

\begin{defi}
  Une perturbation $f$ de $f_0$ {\em relativement aux liens  d'un  îlot stochastique de $f_0$}\footnote{Souvent, nous ne rappellerons pas explicitement l'îlot stochastique concerné et parlerons de perturbation de $f_0$ relativement aux liens.} bordé par $m$ liens hétéroclines est une perturbation telle que les branches stables ou instables qui définissent le bord de l'îlot stochastique de $f_0$ ont leurs prolongements par continuité qui définissent encore $m$ liens hétéroclines  (pour $f$). 
  
  Un îlot stochastique $\Ic_0$ pour $f_0$ est {\em robuste relativement aux liens} s'il existe un voisinage $\Uc$ de $f_0$ pour la topologie $C^2$ tel que tout difféomorphisme (non nécessairement conservatif) de $f\in\Uc$ qui est une perturbation de $f_0$ relativement aux liens de $\Ic_0$ a un îlot stochastique bordée par les  liens perturbés.

\end{defi}

\begin{prop}\label{Probustelien}
L'exemple construit par Berger et Turaev est robuste relativement aux liens.
\end{prop}
Pour démontrer ceci, Berger et Turaev s'inspirent  d'un résultat non publié de Arroyo et Pujals concernant les difféomorphismes des surfaces à bord qui sont   structurellement stables par $C^k$-perturbation\footnote{Un diff\'eomorphisme est  structurellement stable par $C^k$-perturbation si toute dynamique  qui est $C^k$-proche d'un tel difféomorphisme lui est $C^0$-conjuguée.}.

Ils commencent pas faire une perturbation de la perturbation $f$ de $f_0$ de façon à ce que les liens de $f_0$ et sa perturbation que nous noterons encore $f$ soient les mêmes\footnote{Comme les aires des deux îles stochastiques ne sont pas forcément les mêmes, cette perturbation ne préserve pas forcément les aires.}. Ensuite, ils se ramènent à travailler dans le tore en utilisant le revêtement ramifié introduit précédemment, et constatent qu'ils sont en train de perturber en topologie Lipschitz un difféomorphisme Anosov, ce qui donne  un homéomorphisme qui est  encore un difféomorphisme sauf aux quatre singularités du revêtement ramifié  et a encore un exposant de Lyapunov positif partout sauf en ces points. Ils passent ensuite au quotient pour revenir à la sphère puis au disque, et la positivité des exposants de Lyapunov est préservée par ces opérations\footnote{Berger et Turaev montrent aussi que la dynamique perturbée est transitive dans l'îlot stochastique rétabli, en adaptant à leur cas la démonstration valable dans le cas Anosov.}. 
\subsection{Lissage des îlots stochastiques}

Donnons un énoncé qui, même s'il n'est pas donné sous cette forme dans leur texte, est démontré par Berger et Turaev. Cet énoncé leur permettra de lisser les exemples qu'ils construiront qui ne seront pas forcément de classe $C^\infty$ dans un premier temps.
\begin{prop}\label{proplissage}
Soit $r\geq 2$ un entier. Soit $f_0\in {\rm Diff}_\omega^r(\D)$ ayant un îlot stochastique sans coin $\cI_0$ robuste relativement aux liens et $\Uc$ un voisinage de $f_0$ en topologie~$C^r$. Il existe $f\in \Uc$ de classe $C^\infty$ qui est une perturbation de $f_0$ relativement aux liens de $\Ic_0$.  \end{prop}

L'idée est la suivante. Il existe $\hat f\in{\rm
  Diff}^\infty_\omega(\D)\cap\Uc$  qui est proche de $f_0$ en topologie~$C^r$,
ainsi qu'une réunion finie de cercles  $C^\infty$ immergés $\bigcup_{1\leq
  k\leq N}\gamma_k$ qui est proche de $\partial\cI_0$ toujours en
topologie~$C^r$. Or, comme chaque composante connexe de  $\partial \cI_0$ est
invariante par $f_0$, alors chaque $\hat f(\gamma_k)$ est proche en
topologie~$C^r$ de~$\gamma_k$. 
On construit alors un élément~$g$ de~${\rm Diff}^\infty_\omega(\D)$ dans un voisinage en topologie~$C^r$ de~$\mathrm{Id}$  qui envoie chaque $\hat f(\gamma_k)$ sur $\gamma_k$\footnote{Pour faire ceci, il est par exemple possible d'utiliser le théorème du voisinage tubulaire symplectique de \textcite{Weinstein1977} qui permet de travailler dans un anneau $\T\times [-\varepsilon, \varepsilon]$.}. Alors, $f=g\circ \hat f\in \Uc$ est  proche de $f_0$ en topologie $C^r$ et fixe la courbe $\gamma_k$ qui est proche en topologie $C^r$ d'une composante connexe de  $\partial \cI_0$.  La dynamique restreinte~$f_{|\cup\gamma_k}$ est alors proche en topologie $C^r$ de la dynamique $f_{0|\partial\Ic_0}$. Or,  $f_{0}$ restreint à une composante connexe de $\partial \cI_0$ est un difféomorphisme d'un cercle immergé qui a des points périodiques qui sont tous hyperboliques. Une telle dynamique est structurellement stable, ce qui signifie que toute dynamique qui lui est $C^1$-proche lui est $C^0$-conjuguée. D'autre part, les points périodiques contenus dans cette courbe sont toujours hyperboliques et de même type que ceux qui existaient avant la perturbation (attractifs ou répulsifs pour la restriction à la courbe). Ceci implique que $\gamma$ est un lien hétérocline pour $f$, et donc que~$f$ a un îlot stochastique proche de $\cI_0$.

 \subsection{Bonnes cartes pour un bi-lien}\label{secbonnecarte}
 Nous allons maintenant  expliquer une notion due à  Berger et Turaev, celle de bonne carte pour un bi-lien. Dans cette carte, au voisinage de certains arcs des branches du bi-lien, la dynamique se lit comme une translation (les auteurs parlent alors de {\em coordonnées   énergie-temps}) et une certaine transition d'un arc à un autre se  lit comme une application affine hyperbolique. Ces cartes seront utilisées dans la section~\ref{seclien} pour évaluer un certain opérateur de transformation de graphe.

 Étant donné un bi-lien $\cI$ sans coin du disque pour un difféomorphisme $f\in {\rm Diff}_\omega^r(\D)$, il existe toujours un difféomorphisme $\phi\in {\rm Diff}^{r-1}_\omega(\D)$ tel que pour $\hat f=\phi\circ f\circ \phi^{-1}$, le bi-lien  $\hat\cI=\phi(\cI)$ intersecte deux bandes verticales $V_a=[x_a-2\tau, x_a]\times\R$ et $V_b=[x_b, x_b+2\tau]\times \R$ de la manière que nous décrivons ci-dessous. On note $L_a$ et $L_b$ les deux branches du bi-lien.
 \begin{itemize}
 \item $V_a\cap L_a=[x_a-2\tau, x_a]\times \{ y_a\}$ et sur un voisinage $N_a$ de $[x_a-2\tau, x_a]\times \{ y_a\}$, $\hat f(x,y)=(x-\tau, y)$ est une translation horizontale;
 \item $V_b\cap L_a=\emptyset$;
 \item $V_b\cap L_b=[x_b, x_b+2\tau]\times\{y_b^-, y_b^+\}$ avec $y_b^-<y_b^+$;
 \item sur un voisinage $N_b^+$ de $[x_b, x_b+2\tau]\times \{ y_b^+\} $, $\hat f(x, y)=(x+\tau, y)$;
  \item sur un voisinage $N_b^-$ de $[x_b, x_b+2\tau]\times \{ y_b^-\} $, $\hat f(x, y)=(x-\frac{\tau}{2}, y)$;
  \item Il existe $n\geq 1$ tel que $\hat f^n$ envoie $[x_b, x_b+2\tau]\times \{ y_b^+\} \cup f^2([x_b, x_b+2\tau]\times \{ y_b^+\} )$ sur $[x_b, x_b+2\tau]\times \{ y_b^-\} $ et $\forall (x, y)\in N_b^+, f^n(x, y)=\Theta-(\frac{x}{2}, 2y)$ où $\Theta$ est un point fixé de $\R^2$. 
 \end{itemize}
  \begin{figure}[H]
\centering\includegraphics[width=6cm]{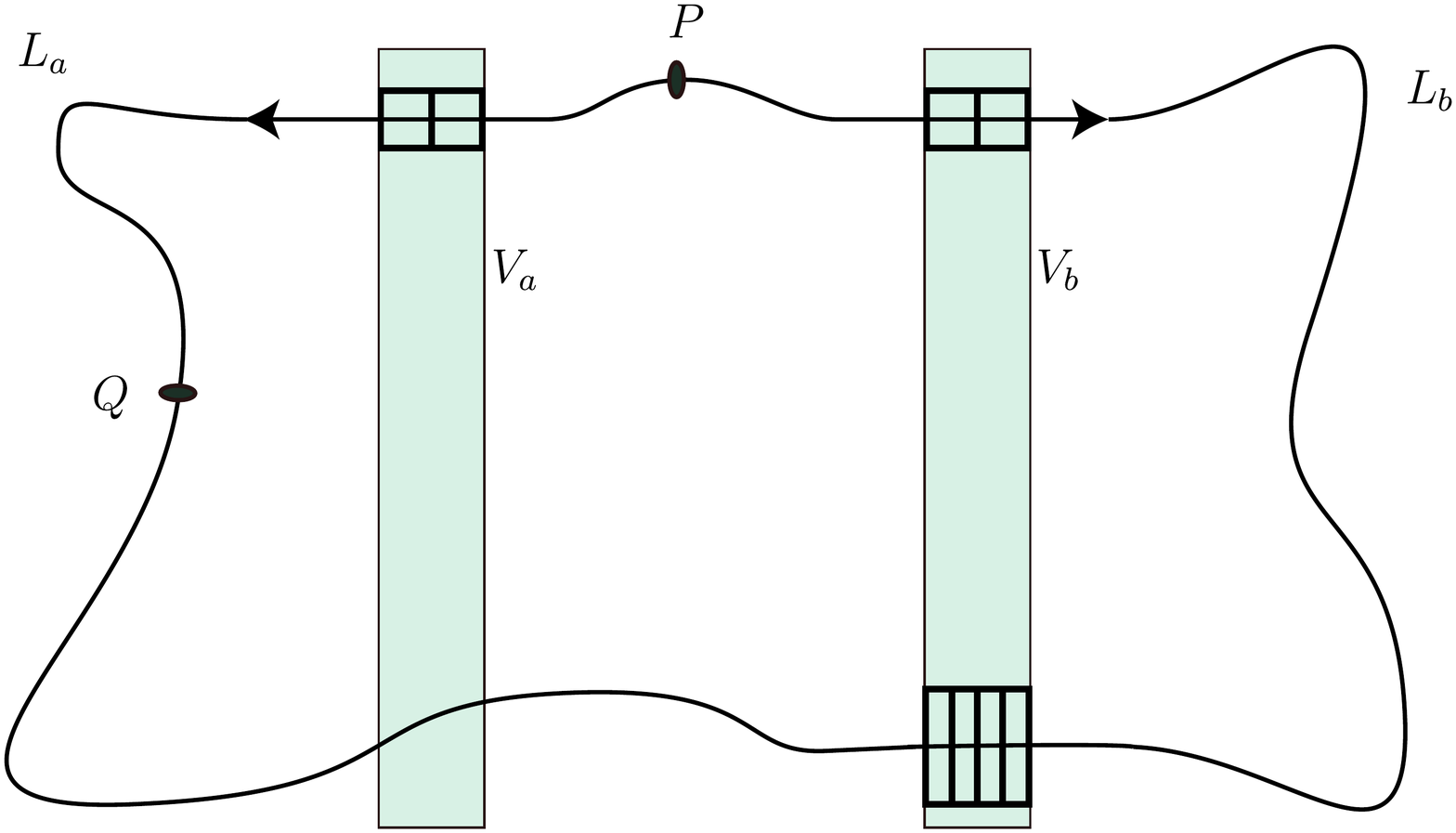}
\caption{Le bi-lien dans une bonne carte.}
\end{figure}

\section{Renormalisation} \label{secrenorm}
\subsection{Première approche}
Expliquons maintenant quel procédé Berger et Turaev vont utiliser pour  injecter leur îlot stochastique (qui encore une fois est associé à  un difféomorphisme  $f_0$ qui est loin de l'identité en topologie $C^1$) dans une dynamique proche de l'identité pour la topologie~$C^\infty$.

Commençons par faire trois remarques qui découlent aisément de la définition de l'entropie métrique.
\begin{itemize}
\item on a $h_\mu(f^n)=nh_\mu(f)$ et donc si une itérée de $f$ est d'entropie métrique positive, il en est de même pour $f$;
\item si $A$ est un ensemble invariant par $f$ tel que $\mu(A)>0$ et si $\mu_A(B)=\frac{\mu(A\cap B)}{\mu (A)}$, alors 
$h_\mu(f)\geq \mu(A).h_{\mu_A}(f_{|A})$ et donc si une restriction de $f$ à un  ensemble invariant est d'entropie métrique positive pour la probabilité conditionnelle, il en est de même pour $f$;
\item l'entropie métrique est invariante par conjugaison qui préserve la mesure\footnote{Kolmogorov l'a d'ailleurs introduite pour montrer que les applications décalages $(u_n)\in \{1, \dots, p\}^\Z\mapsto (u_{n+1})$ munies de la  mesure  produit des  probabilités uniformes ne sont pas mesurablement conjuguées pour des entiers $p$ différents.}.
\end{itemize}
On va donc essayer d'écrire ce difféomorphisme $f_0$ qui est loin de l'identité comme $f_0=h^{-1}\circ g^n_{|A} \circ h$ où
\begin{itemize}
\item $g$ est un difféomorphisme du disque qui préserve l'aire et est proche de ${\rm Id}$ en topologie $C^\infty$;
\item $A$ est un ensemble invariant par $g^n$  d'aire $a>0$;
\item $h$ est un difféomorphisme de classe $C^\infty$ qui multiplie la forme d'aire par $a$.
\end{itemize}
\begin{defi}\label{defrenorm}
Étant donné un difféomorphisme symplectique $f$, on appelle {\em renormalisation} de $f$ tout difféomorphisme qui s'écrit $h^{-1}\circ f^n_{|A} \circ h$ où $h:\D\rightarrow A$ est une application de jacobien constant du disque sur l'ensemble $f^n$-invariant $A$.

Par extension, si $f\in {\rm Diff}_\omega(\D)$, pour $h:U\rightarrow \D$ application de jacobien constant définie sur un ouvert $U$ contenant $\D$, à condition que $f^n(h(\D))\subset h(U)$, on dira que $h^{-1}\circ f^n\circ h_{|\D}$ est une renormalisation de $f$. \end{defi}
Remarquons que comme $h(\D)$ n'est dans ce cas pas forcément invariante par $f^n$ dans la dernière définition, cette notion de renormalisation ne rend pas compte d'une dynamique sous-jacente et ne permet en aucune façon de calculer l'entropie. 

On veut  trouver un difféomorphisme $g$ du disque qui est proche de l'identité en topologie $C^\infty$, préserve l'aire et a un disque topologique $D$  invariant  par un de ses itérés~$g^n$ tel que $g^n_{|D}$ est conformément conjuguée à $f_0$.
\begin{figure}[H]
\centering\includegraphics[width=10cm]{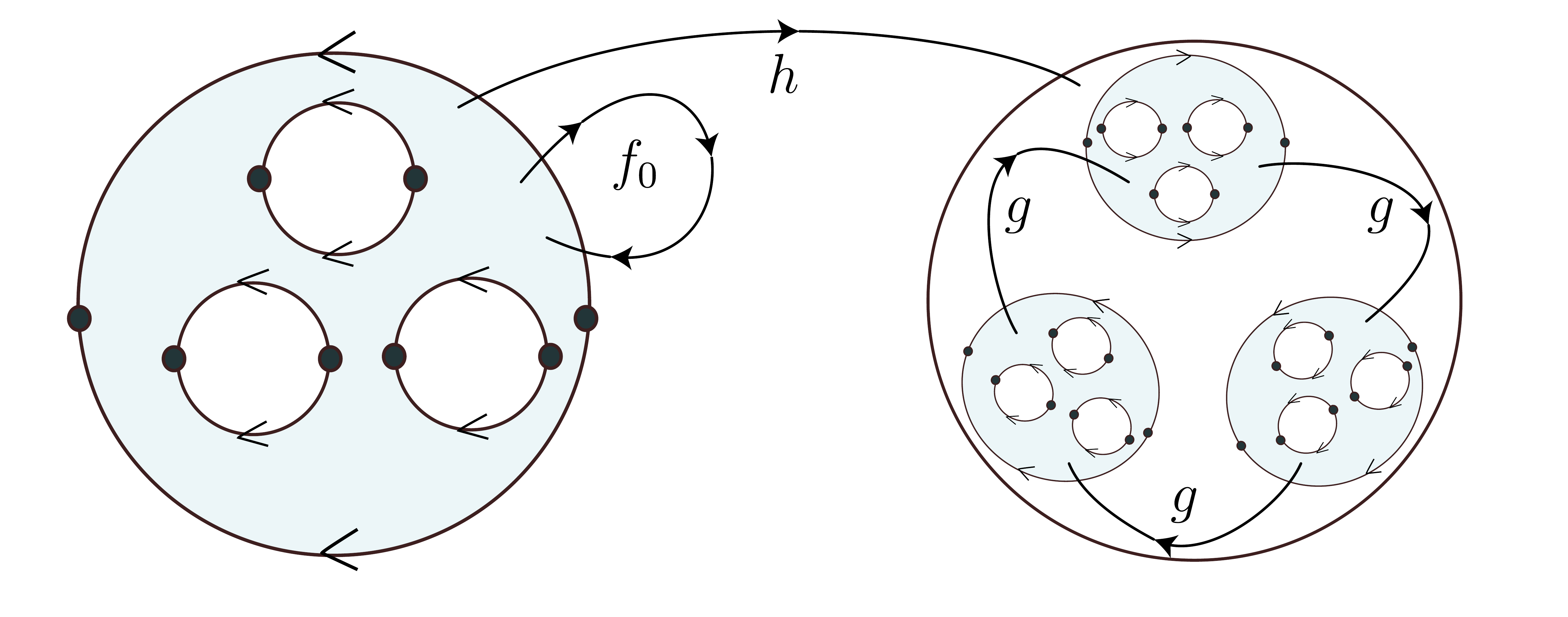}
\caption{Une renormalisation avec $n=3$.}
\end{figure}

Une première approche pour construire $g$ est d'écrire $f_0$ comme la composée $f_0= g_n\circ \dots \circ g_1$ d'un certain nombre de difféomorphismes $g_1, \dots, g_n$  qui sont $\frac{1}{n}$ proches de l'identité en topologie $C^\infty$ et préservant l'aire, ce qui est toujours possible en utilisant une isotopie hamiltonienne.  On place dans le disque $\D$ $n$ petits disques disjoints $D_1, \dots, D_n=D_0$ de même rayon $r$ (qui est donc plus petit que $\frac{1}{\sqrt{n}}$ car la somme de leurs aires est plus petite que celle du disque $\D$). Notons $h_j:\D\rightarrow D_i$ l'homothétie de rapport $r$ qui envoie $\D$ sur $D_j$. On construit  alors un difféomorphisme  $f$ de $\D$ tel que $f(D_j)=D_{j+1}$ et $f_{|D_j}=h_{j+1}\circ g_j\circ h_j^{-1}$.  Par construction, on a $h_1^{-1}\circ g^n_{|D_1}\circ h_1=f_0$ et donc $g$ est bien d'entropie métrique positive. Comme les $h_j$ sont des homothéties de rapport $r$, la distance en topologie $C^k$ de  $f_{|D_j}=h_{j+1}\circ g_j\circ h_j^{-1}$ à l'identité  est de l'ordre de $\frac{1}{r^{k-1}}\| g_j\|_{C^k}\geq (\sqrt{n})^{k-1}\| g_j\|_{C^k}$.

Le théorème des accroissements finis nous dit alors que $\| g_i\|_{C^1}$ est au moins de l'ordre de $\frac{1}{n}$, donc $\| f_{|D_j}\|_{C^k}$ est au moins de l'ordre de $(\sqrt{n})^{k-3}$, donc explose dès que $k>3$.

Si nous n'étions intéressés qu'en des approximations en topologie $C^1$ ou $C^2$, ce raisonnement serait suffisant.  C'est ce que font Newhouse, Ruelle et Takens pour construire des difféomorphismes du tore qui ont des attracteurs étranges dans des voisinages petits en topologie $C^2$ de l'identité, \parencite{NewhouseRuelleTakens1978}. Mais nous sommes intéressés en des approximations plus fines.  

Maintenant, nous allons expliquer comment la renormalisation au voisinage des bandes homoclines permet à Berger et Turaev de réaliser beaucoup de dynamiques dans des topologies plus fines.  

\subsection{Énoncé du lemme de renormalisation au voisinage des bandes homoclines}
Soit $(M, \omega)$ une surface munie d'une forme d'aire et $f:M\rightarrow M$ un difféomorphisme symplectique de classe $C^\infty$ qui a un point périodique hyperbolique $O$ de période $\tau$. Cela signifie que  $f^\tau (O)=O$ et les deux valeurs propres $\lambda\in ]0, 1[$, $\lambda^{-1}$ de $Df^\tau(O)$ sont réelles et de module différent de $1$.  On sait alors qu'il existe deux immersions injectives $j_s, j_u:(\R, 0)\hookrightarrow (M, O)$ de classe $C^\infty$ telles que $j_s(\R)$ (resp. $j_u(\R)$) est la variété stable (resp. instable) de $O$ qui est définie par
$$W^s(O, f)=\{ x\in M; \lim_{n\rightarrow +\infty} d(f^nx, f^nO)=0\}=W^u(O, f^{-1}).$$
Les deux variétés immergées $W^s(O, f)$ et $W^u(O, f)$ sont invariantes par $f^\tau$ et se coupent transversalement en $O$. Parfois,  $W^u(O, f)$ et $f^j(W^s(O, f))=W^s(f^j(O), f)$ se rencontrent en un autre point que $O$ pour un  $j\in [0, \tau-1]$. On dit alors que ces points d'intersection sont des points homoclines à l'orbite de $O$~: l'orbite positive et l'orbite négative d'un tel point se rapprochent de celle de $O$ quand le temps tend vers l'infini. Il peut arriver qu'en  un tel point les deux variétés se coupent transversalement, et dans ce cas on a un fer à cheval. À l'autre extrême, ces deux variétés peuvent être égales ou au moins coïncider le long d'un arc $I$ non réduit à un point. On parle alors de {\em bande homocline}. 
\begin{figure}[H]
\centering\includegraphics[width=6cm]{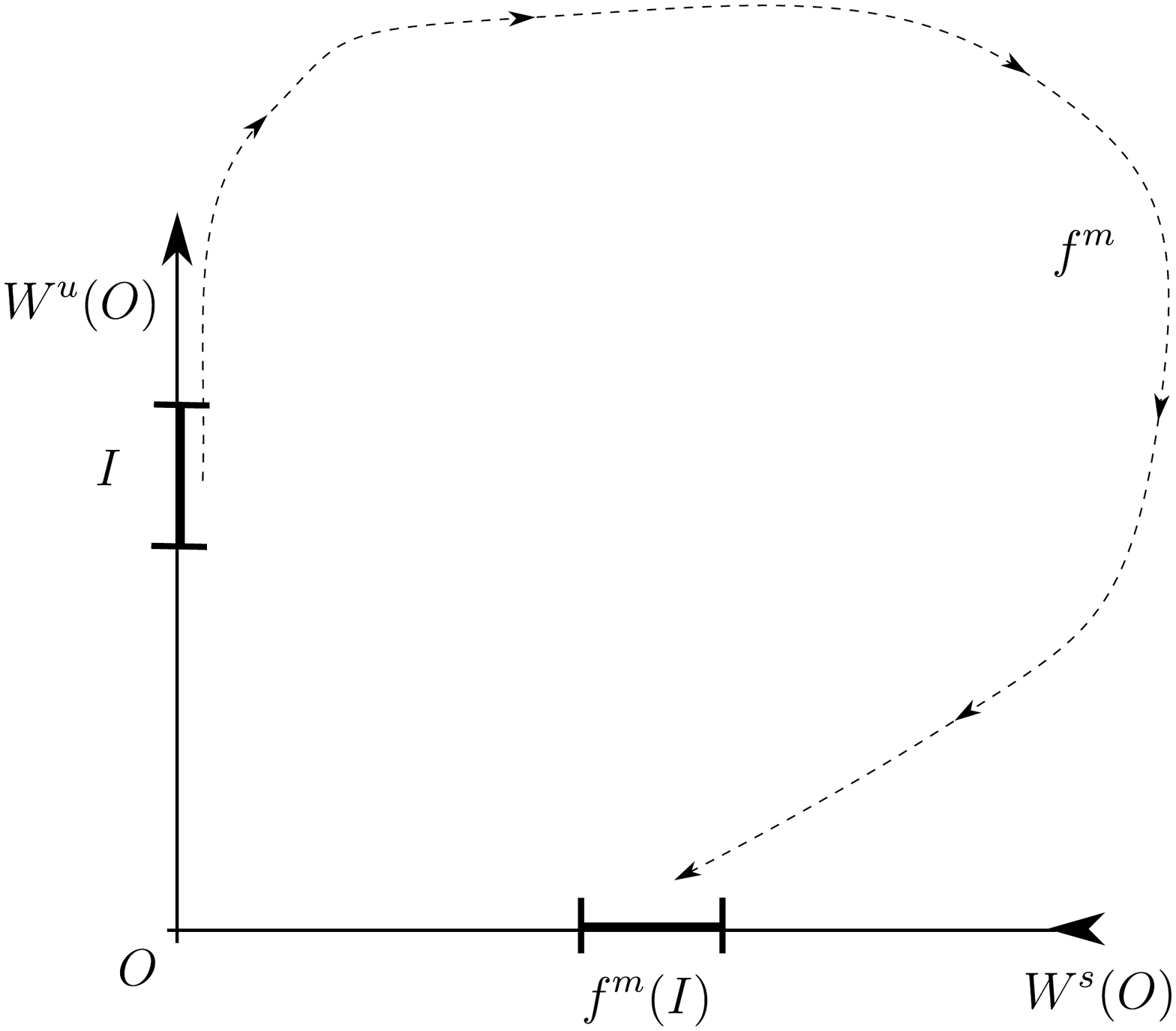}
\end{figure}
Nous supposerons désormais l'existence d'une telle bande homocline. En se plaçant dans une bonne carte, on peut supposer que  les variétés stables et instables de $O$ coïncident localement avec les axes des coordonnées et donc que le segment vertical $I$ s'envoie par une itérée $f^m$ de $f$ sur un segment horizontal.  Sous ces hypothèses, Berger et Turaev développent \textcite{GonchenkoTuraevShilnikov2007} et démontrent le lemme suivant de renormalisation.

\begin{defi}
Soit $\psi:\R\rightarrow \R$ une application de classe $C^r$. L'application de type Hénon $H_\psi$ est définie par
$$H_\psi(x, y)=(y, -x+\psi(y)).$$
\end{defi}
\begin{lemm}[Lemme de renormalisation]\label{lemrenorm} Pour tout $N$ impair, pour tout $r\geq 1$ pour tous $L, \delta>0$ et pour tout voisinage $U_I$ de $I$, il existe
\begin{itemize}
\item un entier $n>1$;
\item des difféomorphismes $\Phi_1, \dots, \Phi_N\in{\rm Diff}^\infty _\omega(\R^2)$ qui sont $\delta$-proches de l'identité;
\item un difféomorphisme de classe $C^\infty$ $Q:\R^2\rightarrow \R^2$ de jacobien constant tel que $Q(\D)\subset U_I$;
\end{itemize}
tels que pour toutes fonctions $\psi_1, \dots , \psi_N\in C^r(\R, \R)$ dont les normes~$C^r$ sont majorées par $L$, il existe un difféomorphisme symplectique $\hat f$ de classe $C^r$\footnote{Si les $\psi_i$ sont de classe $C^\infty$, $\hat f$ est de classe $C^\infty$ aussi mais on n'a un contrôle que sur la norme~$C^r$.} qui coïncide avec f sur le complémentaire de $U_I$, est tel que $\| \hat f-f\|_{C^r}<\delta$ et 
$$Q^{-1}\circ \hat f^n\circ Q_{|\D}=H_{\psi_N}\circ \Phi_N\circ \dots \circ H_{\psi_1}\circ \Phi_{1|\D}.$$
De plus, $\hat f^m\circ Q(\D)\cap Q(\D)=\emptyset$ pour tout $m=1, \dots , n-1$.
\end{lemm}
En d'autres termes, il existe $N$ applications $\Phi_1, \dots, \Phi_N$ $C^r$-proches de l'identité telles que, quitte à s'autoriser une perturbation petite $\hat f$ de $f$, on peut réaliser   comme renormalisation de $\hat f$\footnote{Renormalisation s'entend ici au sens faible de la définition~\ref{defrenorm}, c-à-d sans que le domaine où on renormalise soit invariant par une itérée de $f$.}  toutes les composées de $N$ applications de Hénon de norme $C^r$ bornée par $L$ alternées avec les $\Phi_i$.

Un résultat de \textcite{Turaev2003} permet d'approximer toute dynamique symplectique par une composée d'applications de Hénon.  Le résultat précédent donne donc une clé pour voir apparaître un ensemble dense de dynamiques symplectiques comme renormalisées.  Dans un autre cadre, 
\textcite{Turaev2003}
introduisait  la notion d'application universelle telle que l'ensemble de ses
applications renormalisées est dense dans ${\rm Diff}^\infty_\omega(\D,
\R^2)$\footnote{Pour des difféomorphismes de classe $C^1$, Bonatti et Diaz
  montraient l'existence de telles dynamiques en présence de classes
  homoclines d'un certains type \parencite{BonattiDiaz2002}.}.
 
 \subsection{Démonstration du lemme de renormalisation}\label{ssdemlemrenorm}
 Quitte à diminuer l'intervalle $I$, on peut supposer que $I\cap f^\tau(I)=\emptyset$ et  que donc $I$ ne rencontre aucune de ses images. Quitte à les remplacer par des images par des itérés de la dynamique, on peut aussi supposer que $I$ et $f^m(I)$ sont dans un petit voisinage $V=[-\alpha, \alpha]^2$ de $O$. On fixe deux voisinages ouverts $V_\varepsilon\subset \bar V_\varepsilon\subset U_\varepsilon$ de taille $\varepsilon>0$ de $f^m(I)$. On choisit alors $N$ points distincts $M_1^- \dots, M_N^-$ dans $I$ et on utilise la notation $M_i^+=f^m(M_i^-)$. On note leurs coordonnées $M_i^-(0, y_i^-)$ et $M_i^+(x_i^+, 0)$.
 
 On va utiliser tout d'abord deux types de changement d'échelle qui sont des affinités\footnote{Rappelons que $\lambda\in]0, 1[$ est une des deux valeurs propres de $Df^\tau(O)$, l'autre étant $\lambda^{-1}$.}
 \begin{itemize}
 \item au voisinage de $M_i^-$, $(x,y)\mapsto (\frac{x}{\lambda^k}, y-y_i^-)$;
 \item au voisinage de $M_i^+$, $(x, y)\mapsto (x-x_i^+, \frac{y}{\lambda^k})$.
 \end{itemize}
 On utilisera aussi une suite $(\mu_k)$ qui tend en décroissant de façon sous-exponentielle vers~$0$ pour faire  des  changements d'échelle homothétiques $(x,y)\mapsto (\frac{x-x_i}{\mu_k}, \frac{y-y_i}{\mu_k})$  autour de points variables $M_i(x_i, y_i)$.
 
 Pour faire ces  changements d'échelle, on a donc besoin d'être 
 \begin{itemize}
 \item soit dans une bande verticale de largeur $\lambda^k\alpha$ et dans une ellipse verticale qui est proche d'un $M^-_i$ de grand axe  au plus $\mu_k\alpha$ et dont le rapport des axes vaut $\lambda^k$;
 \item soit dans une bande horizontale de hauteur $\lambda^k\alpha$ et dans une ellipse horizontale qui est proche  d'un $M^+_i$ de grand axe au plus $\mu_k\alpha$ et dont le rapport des axes vaut $\lambda^k$\footnote{Berger et Turaev font  aussi un  changement d'échelle indépendant de $k$ que nous ne faisons que mentionner pour ne pas alourdir.}.
 \end{itemize}
  Un premier résultat concerne l'application $f^{\tau k}$ renormalisée dans une ellipse de taille au plus $\mu_k\alpha$ dont les $k$ premières images par $f^\tau$ restent dans $V$ \parencite{GonchenkoShilnikovTuraev2008}.
   \begin{figure}[ht]
\centering \includegraphics[width=4cm]{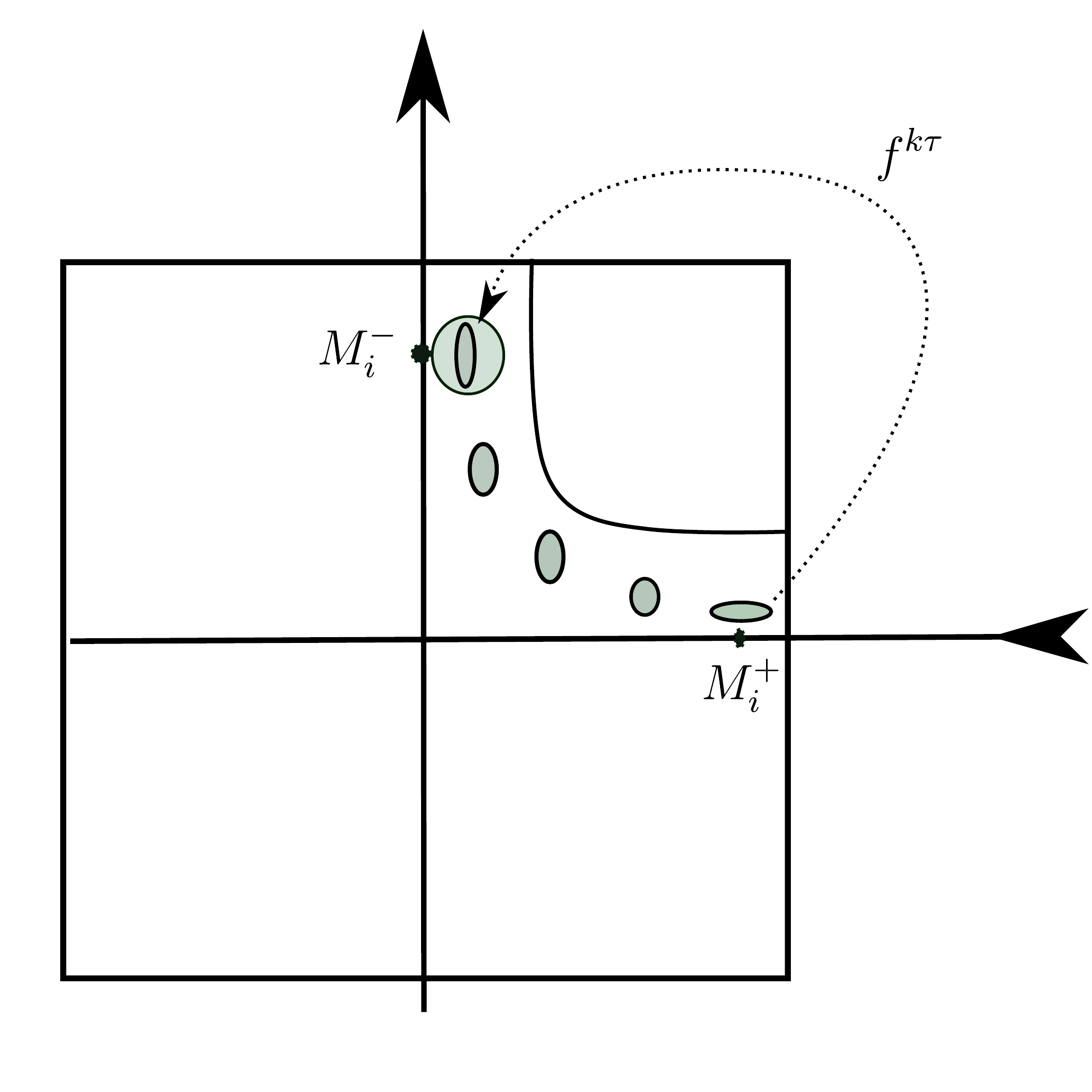}
\caption{Une transition \og presque linéaire\fg}
\end{figure}
   Il énonce que l'application $f^{\tau k}$ renormalisée dans une ellipse  de taille au plus $\mu_k\alpha$ dont les $k$ premières images par $f^\tau$ restent dans $V$ s'écrit (dans les coordonnées renormalisées  fixées de taille $R$)
 $$T_0(X, Y)=(X, Y)+\Delta_k(X,Y)$$
 où la norme $C^\infty$ de $\Delta_k$ tend vers $0$ quand $k$ tend vers
 $+\infty$. En d'autres termes, les transitions dans $V$, où $f^\tau$ est \og
 presque linéaire\fg, se lisent après renormalisation comme \og presque l'identité\fg.

 Un second résultat  décrit ce qui se passe pour l'application renormalisée quand on suit une connexion homocline plate de $M_i^-$ à $M_i^+$ par $f^m$. 
\begin{figure}[H]
\centering \includegraphics[width=4cm]{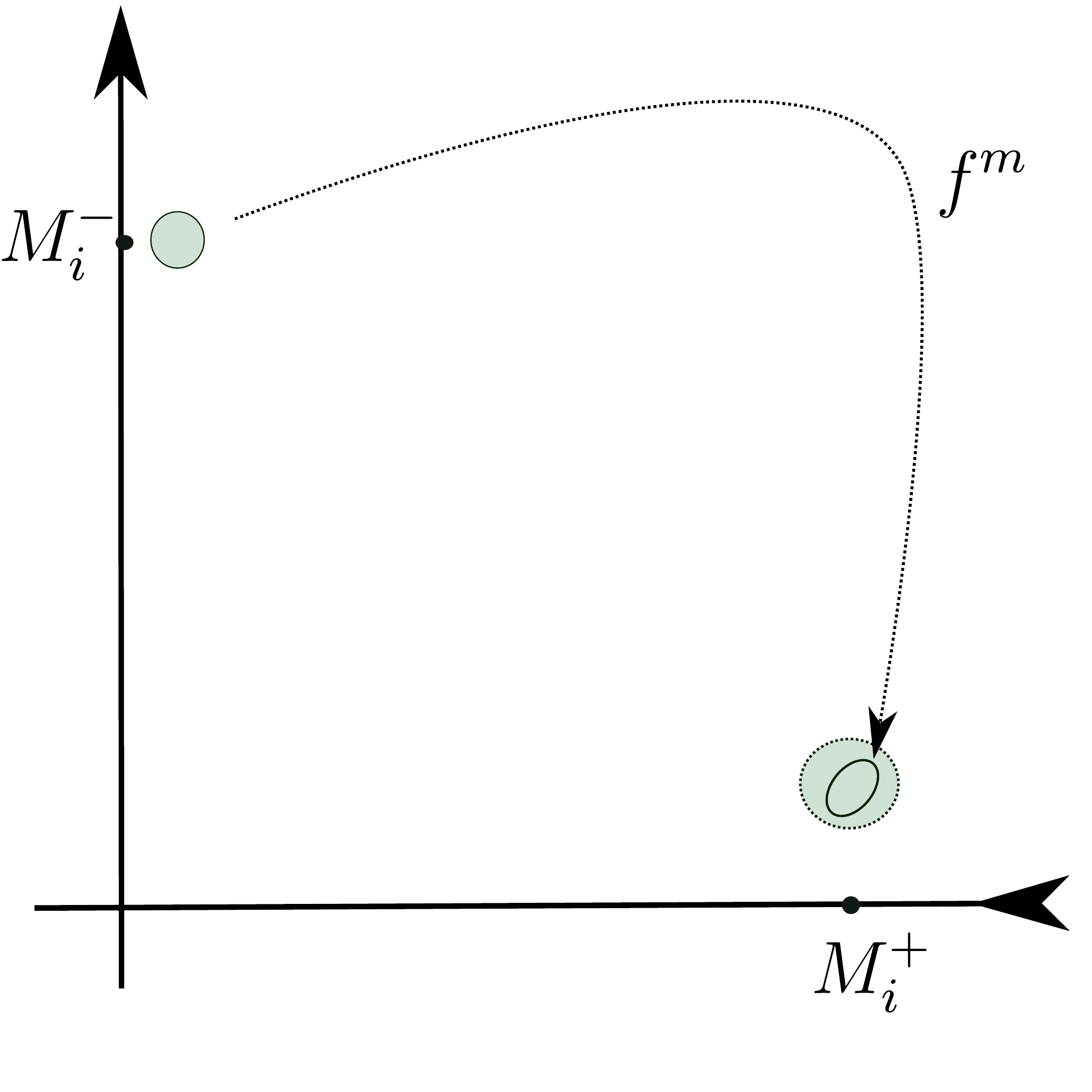}
\caption{Une autre transition \og presque linéaire\fg}
\end{figure}

En ce cas, l'application renormalisée n'est plus proche de l'identité, mais d'une application de Hénon linéaire, toujours pour la topologie~$C^\infty$.
$$(X, Y)\mapsto (Y, -X+aY);$$
avec $a$ uniformément borné en $k$.

Berger et Turaev vont alors faire une perturbation à  support très petit de $f$ de manière à intercaler entre des transitions bien choisies des applications de la forme $(x, y)\mapsto (x, y+\psi(x))$\footnote{Il s'agit  donc de translations dans la fibre verticale.}qui permettront de faire apparaître des applications de type Hénon et aussi de faire disparaître la transition qui n'est a priori pas proche de l'identité. Pour ce faire, il faut préciser un peu les domaines des transitions et des renormalisations.

Tout d'abord, ils choisissent   des ellipses horizontales $B_{i, k}^+$ au voisinage des $M_i^+$ et des ellipses verticales $B_{i, k}^-$ au voisinage des $M_i^-$ de telle sorte que $f^{\tau k}$ définisse une transition dans $V$ de $B_{i, k}^+$ vers $B_{i+1, k}^-$\footnote{Notons que de telles ellipses de transition peuvent être construites quels que soient les points $(M_i^+)$ et ($M_i^-$) fixés.}. Le centre de $B_{i, k}^+$ a même abscisse $x_i^+$ que $M_i^+$ et le centre de $B_{i+1, k}^-$ a même ordonnée que $M_{i+1}^-$. Si la transition envoie le centre de  $B_{i, k}^+$  sur le centre de $B_{i+1, k}^-$, ceci détermine entièrement les centres des ellipses.
\begin{figure}[H]
\centering \includegraphics[width=5cm]{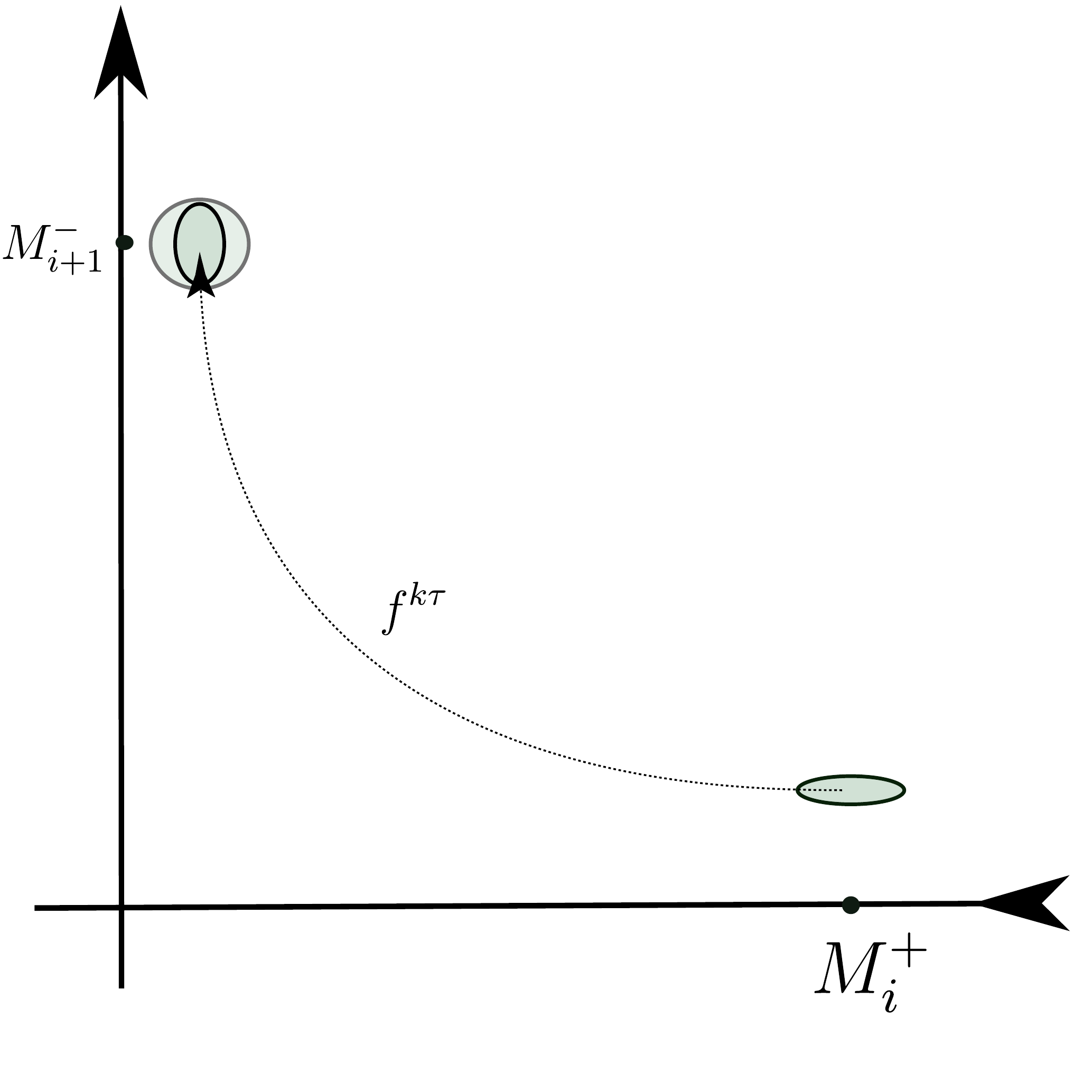}
\caption{Transition du voisinage de $M_i^+$ vers le voisinage de $M_{i+1}^-$.}
\end{figure}
On remarque que $f^{\tau k}(B_{i, k}^+)$ peut avoir un diamètre plus grand que celui de $B_{i, k}^+$, c'est pourquoi le diamètre  de $B_{i+1, k}^-$ peut être plus grand que celui de $B_{i, k}^+$. Mais Berger et Turaev montrent que le rapport entre les deux rayons reste borné uniformément en $k$.

Ensuite, on utilise $f^m$ pour faire une transition de $B_{i+1, k}^-$ vers un voisinage de $M_{i+1}^+$.  Il n'y a aucune raison que $f^m(B_{i+1, k}^-)$ soit incluse dans $B_{i+1, k}^+$, non pas pour une question de taille car celle-ci peut être ajustée\footnote{Là encore le rapport des deux rayons peut être borné uniformément  en $k$.}, mais parce que le centre de $B_{i+1, k}^+$ peut être relativement éloigné de $f^m(B_{i+1, k}^-)$, où le terme {\em relativement} concerne ce qui se passe quand on a renormalisé à la fois par $\lambda^k$ et par $\mu_k$. Toutefois, si on revient dans le voisinage $V$ de $O$, avant les changements d'échelle  donc, on constate que cette distance ne peut être d'ordre plus grand que $\lambda^k$. 

Berger et Turaev construisent alors un difféomorphisme symplectique $g_k$ qui est $\mu_k$-proche de l'identité en topologie~$C^k$, à support dans $U_\varepsilon$, qui coïncide avec une application $(x, y)\rightarrow (x, y+\beta_k(x))$ dans $V_\varepsilon$, et qui a les effets suivants\footnote{Faire ces opérations ne requiert aucune condition spéciale sur les points $M_i^-$ à part qu'ils soient deux à deux distincts.}:
\begin{itemize}
\item pour chaque $i\in [1, N-1]$, on a $g_k\circ f^m(B_{i, k}^-)\subset B_{i+1, k}^+$; on peut en effet grâce  à une translation dans la fibre verticale ramener $f^m(B_{i, k}^-)$ dans $B_{i+1, k}^+$;
\item après renormalisation par $\lambda^k$ et par $\mu_k$, l'application $g_k\circ f^m$ restreinte à $B_{i, k}^-$ est proche de $H_{\psi_i}:(X, Y)\mapsto (Y, -X +\psi_i(Y))$; en effet, toujours grâce à une translation dans la fibre, on peut faire disparaître le terme $a$ qui apparaissait dans la transition homocline $f^m$; on est alors ramené à composer une application $\Phi_i$ proche de $H_0:(X, Y)\mapsto (Y, -X)$ avec l'application $S_{\psi_i}~: (X, Y)\mapsto (X, Y+\psi_i(X))$, ce qui donne la composée d'une application proche de l'identité\footnote{ Notons que cette application $\Phi_i$ dépend de $k$ mais pas de $\psi_i$} et d'une application proche de $(X, Y)\mapsto (Y, -X+\psi_i(Y))$.
\end{itemize} 
La renormalisée d'une composition étant la composition des renormalisées, on peut alors exprimer la renormalisée de $\left((g\circ f)^{N(m+\tau k)}\right)_{|B_{1, k}^+}$ qui est 
$$H_{\psi_N}\circ \hat\Phi_N\circ \dots \circ H_{\psi_1}\circ \hat\Phi_{1|\D}.\footnote{Dans cette expression, les applications $\hat\Phi_i$ ne dépendent pas de $\psi_i$.}$$
Notons que l'image de cette application n'a aucune raison d'être dans $B_{1, k}^+$, on sait juste qu'elle est dans $f^m(B_{1, k}^-)$. Ce n'est donc a priori pas une application de premier retour dans $B_{1, k}^+$.  
\begin{figure}[H]
\centering \includegraphics[width=6cm]{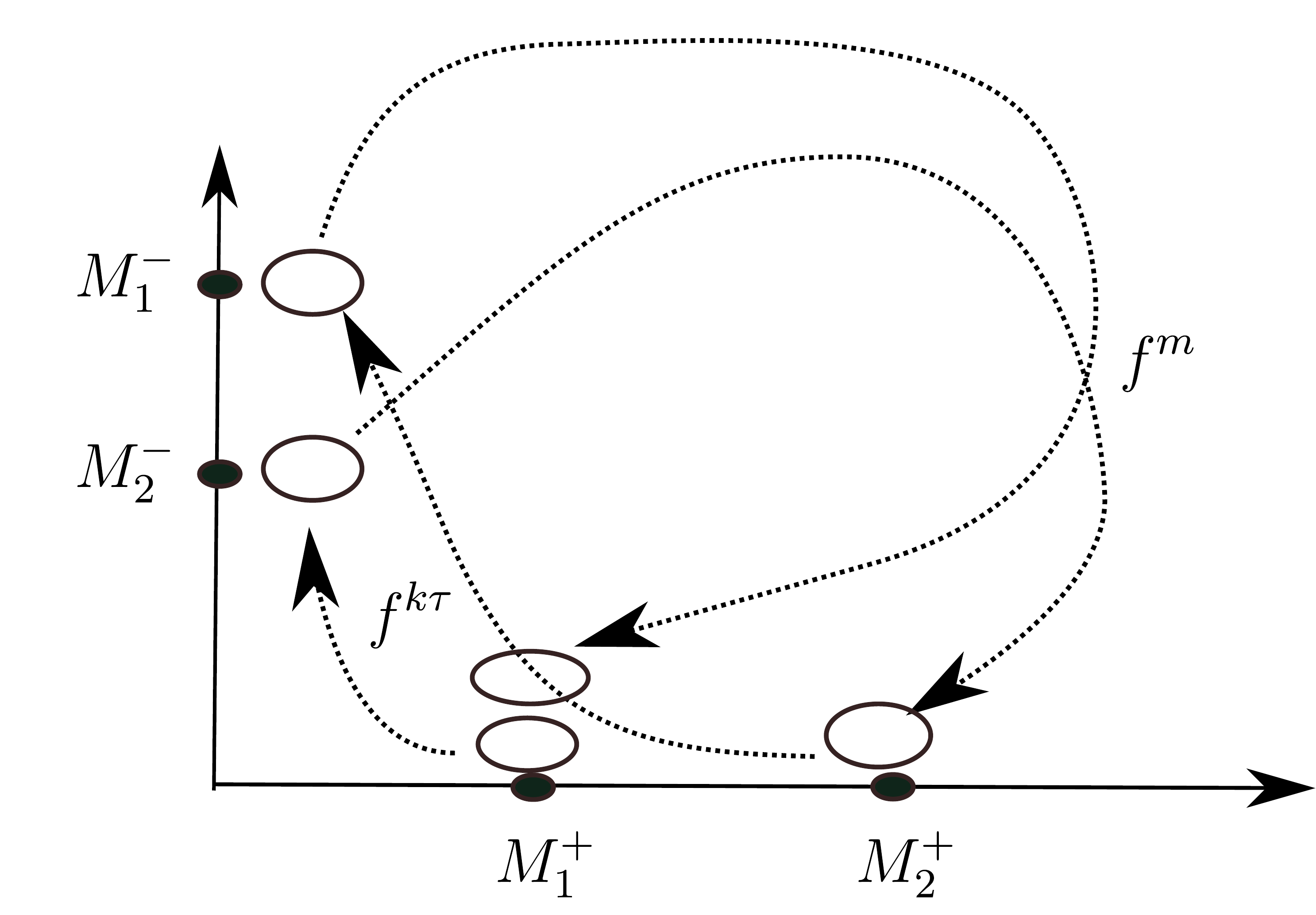}
\caption{la construction pour $N=2$.}\label{figtransitionglobale}
\end{figure}
\subsection{Le corollaire du lemme de renormalisation}
 Berger et Turaev vont en fait utiliser le corollaire suivant du lemme de renormalisation.
\begin{coro}\label{cororenorm}
Soit $ f \in {\rm Diff}_\omega^\infty(M)$ qui admet un point périodique hyperbolique ayant une bande homocline. Alors il existe un sous-ensemble $C^\infty$-dense $\cF$ de ${\rm Diff}_\omega^\infty(\D, \R^2)$ tel que pour tout $F\in \cF$, tout $r\geq 2$, pour toute fonction $\psi: \R\rightarrow \R$  de classe $C^r$ et tout $\varepsilon>0$, il existe un difféomorphisme  $ \hat f \in {\rm Diff}_\omega^ r(M)$ tel que 
\begin{itemize}
\item $d_{C^r} (\hat f, f)<\varepsilon$;
\item la composée  $S_\psi\circ F$ est égale à une   application  renormalisée de $\hat f$ où $S_\psi(x, y)=(x, y+\psi(x))$.
\end{itemize}
\end{coro}
 
Pour déduire le corollaire~\ref{cororenorm} du lemme de renormalisation, les deux arguments clés sont
\begin{itemize}
\item le résultat de \textcite{Turaev2003} qui énonce qu'on peut approximer en topologie~$C^\infty$ tout plongement symplectique du disque dans $\R^2$ par une composée d'un nombre pair d'applications de Hénon polynomiales;
\item la remarque que $S_\psi=H_\psi\circ H_0^{-1}$.
\end{itemize}
Étant donné $F\in {\rm Diff}_\omega^\infty(\D, \R^2)$,  ils  approximent  $H_0^{-1}\circ F$ par une composition $H_{\psi_{2k}}\circ\dots \circ H_{\psi_1} $ d'applications de Hénon polynomiales avec une précision de $\frac{\varepsilon}{2}$ en topologie~$C^\infty$. Comme $H_0(x,y)=(y, -x)$ est une rotation,  on a aussi
$$d_{C^\infty}(F, H_0\circ H_{\psi_{2k}}\circ \dots \circ H_{\psi_1})<\frac{\varepsilon}{2}.$$
Utilisant   le lemme de renormalisation pour $N=2k+1$, un majorant~$L$ de la norme~$C^r$ de $\psi_{2k+1}=\psi, \psi_{2k},  \dots,\psi_1$, $\delta>0$, on  déduit l'existence d'un entier $n\geq 1$, d'un difféomorphisme de classe $C^r$ $\hat f$ qui coïncide avec $f$ en dehors de $U_I$, tel que $d_{C^r}(f, \hat f)<\varepsilon$, de difféomorphismes $\Phi_1, \dots , \Phi_N\in {\rm Diff}^\infty_\omega(\R^2)$ qui sont $\delta$-$C^r$-proches de ${\rm Id}$, un  plongement $Q:\D\rightarrow \R^2$ de classe $C^\infty$ et de jacobien constant tel que 
$$Q^{-1}\circ \hat f^n\circ Q_{|\D}=H_{\psi_N}\circ \Phi_N\circ \dots \circ H_{\psi_1}\circ \Phi_{1|\D}.$$
Comme $H_\psi =S_\psi\circ H_0$, cela s'écrit aussi
$$Q^{-1}\circ \hat f^n\circ Q_{|\D}=S_\psi\circ H_0\circ \Phi_N\circ \dots \circ H_{\psi_1}\circ \Phi_{1|\D}.$$
Posant $\hat F= H_0\circ \Phi_N\circ \dots \circ H_{\psi_1}\circ \Phi_1$ et choisissant $\delta$ petit, on a alors
$$d_{C^\infty}(\hat F, H_0\circ H_{\psi_{2k}}\circ \dots \circ H_{\psi_1})<\frac{\varepsilon}{2}$$
donc $d_{C^\infty}(F, \hat F)<\varepsilon$, $d_{C^r}(f, \hat f)<\varepsilon$ et 
 $$Q^{-1}\circ \hat f^n\circ Q_{|\D}=S_\psi\circ \hat F.$$

\section{La  réparation de liens}\label{seclien}
\begin{figure}[H]
\centering \includegraphics[width=5cm]{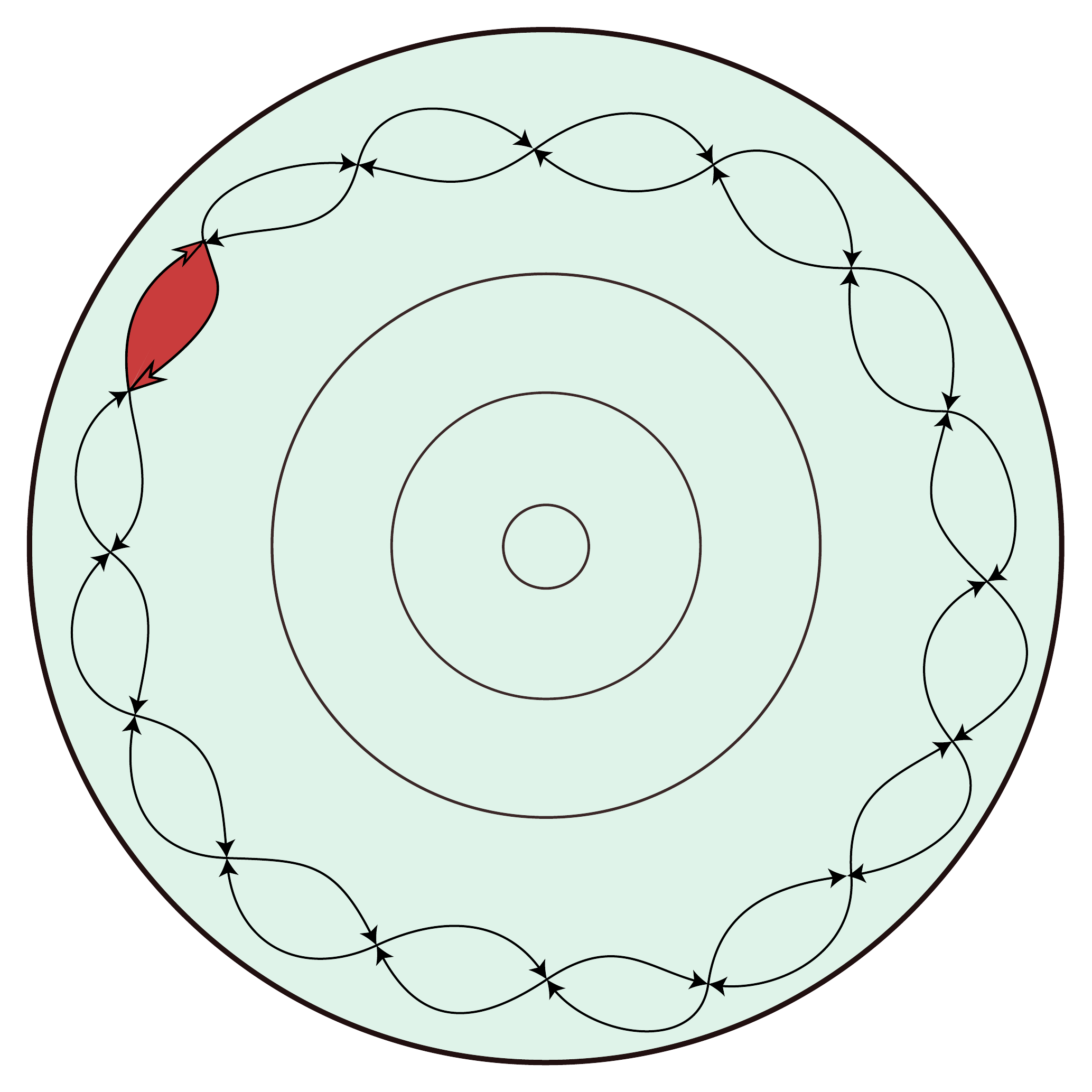}
\caption{Une guirlande de connexions hétéroclines.}\label{figguirlande}
\end{figure}

Afin d'utiliser les résultats obtenus en section~\ref{secrenorm} au voisinage des bandes homoclines, on perturbe l'application identité du disque en une application $f_0$ qui préserve l'aire, vaut l'identité au bord du disque et qui a un point hyperbolique tel que la variété stable et la variété instable de son orbite coïncident\footnote{En utilisant les formes normales,  
  \textcite{GelfreichTuraev2010} ont montré qu'une telle perturbation est possible au voisinage de tout point fixe elliptique. Ceci explique comment montrer le remords cité dans l'introduction.}. Ceci ne présente pas de difficulté majeure~: en utilisant un hamiltonien du type énergie cinétique plus énergie potentielle dans un petit anneau, on fait apparaître des suites de connexions hétéroclines avec des branches stables et instables qui coïncident comme sur la figure \ref{figguirlande} et on compose avec une rotation d'angle petit.

\begin{figure}[H]
\centering \includegraphics[width=8cm]{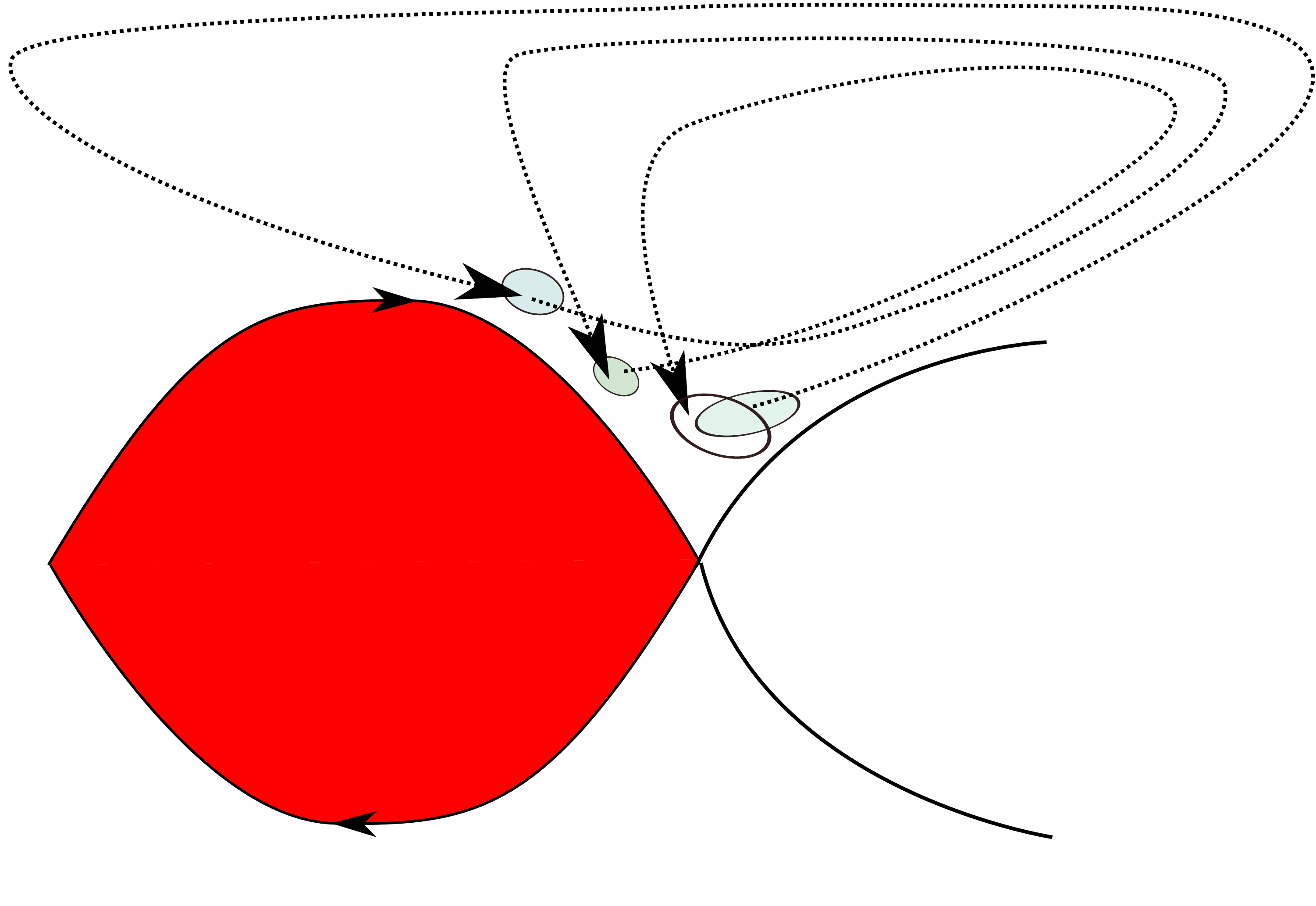}
\caption{La dynamique perturbée itérée sur un petit domaine.}\label{figrenormalisation11}

\end{figure}

Le corollaire~\ref{cororenorm} nous permet alors grâce à une  perturbation $ f
$ de la dynamique $f_0$ de voir une perturbation de l'îlot stochastique $\hat
\cI_0$ qui est l'îlot  introduit en section~\ref{secile}  lu dans une bonne
carte que nous allons décrire un peu plus loin comme une application
renormalisée de $ f$. Nous représentons ci-après ce qui se passe pour une
itérée de $f$  pr\`es du domaine invariant par cette itérée que nous avons
représenté en foncé sur la figure \ref{figguirlande}. La figure \ref{figrenormalisation11} reprend la 
figure \ref{figtransitionglobale} adaptée à  notre cas.

\begin{figure}[H]
\centering \includegraphics[width=6cm]{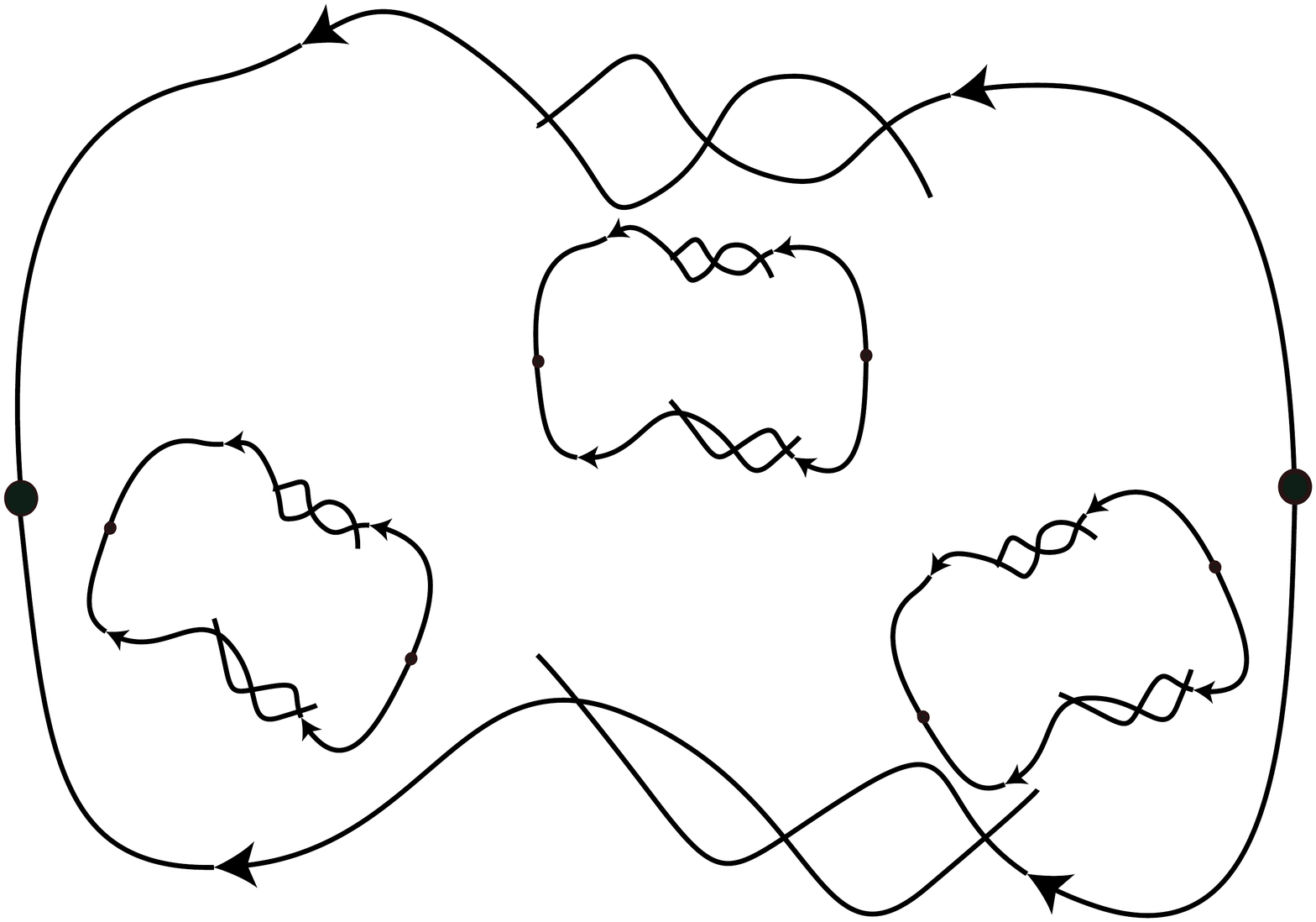}
\caption{L'îlot stochastique perturbé\protect\footnotemark.}\label{figliencasse}
\end{figure}

Le défi est maintenant de restaurer les liens.  Les quatre bi-liens ont en effet été cassés et la figure renormalisée que l'on obtient ressemble à la figure \ref{figliencasse}.
\footnotetext{Le dessin montre des branches stables et instables qui continuent à s'intersecter, il n'y a a priori aucune raison que ce soit le cas et le dessin pourrait être différent.}
Il est assez fréquent en dynamique de   vouloir créer des intersections
 hétéroclines dans le but d'obtenir des fers à cheval et donc un comportement
 chaotique, par exemple quand une branche stable d'un point périodique
 hyperbolique s'accumule sur une branche instable d'un autre point
 périodique. C'est un problème perturbatif difficile dont une solution
 générale existe uniquement en topologie $C^1$, présentée par Hayashi lors du
 congrès international de 1998 \parencite{Hayashi1998}. 
 
 Le problème envisagé ici est différent, car on sait que la dynamique est proche de celle d'un îlot stochastique où branches stables et instables coïncidaient. Le problème est de rétablir cette coïncidence en n'utilisant qu'un certain type de perturbations autorisées par le corollaire~\ref{cororenorm}, celles de la forme $(x, y)\mapsto (x, y+\psi(x))$.\footnote{Remarquons à l'opposé que pour une dynamique holomorphe et entière, 
   \textcite{Ushiki1980} a montré qu'il n'est pas possible d'avoir des bi-liens.
 }
 \subsection{Cartes énergie-temps pour des perturbations de bi-liens}\label{secenergietemps}
 Afin d'estimer la non-coïncidence des variétés stables et instables, Berger et Turaev vont encore  utiliser des coordonnées de type énergie-temps comme en section~\ref{secbonnecarte} dans lesquelles les branches stables et instables s'écriront localement comme des graphes.  Notons que si la dynamique est de classe $C^r$, les cartes énergie-temps construites par Berger et Turaev sont de classe $C^{r-1}$. Alors, les graphes coïncident si et seulement si la différence des fonctions qui leur sont associées est nulle. 
 
On  se ramène donc  à travailler dans des sous-espaces fixés  d'espaces de fonctions  de la forme $C^{r-1}([x, x+\tau], \R)$, deux tels espaces fonctionnels pour chacun  des bi-liens  bordant l'îlot stochastique initial.
 
 \begin{figure}[H]
\centering \includegraphics[width=8cm]{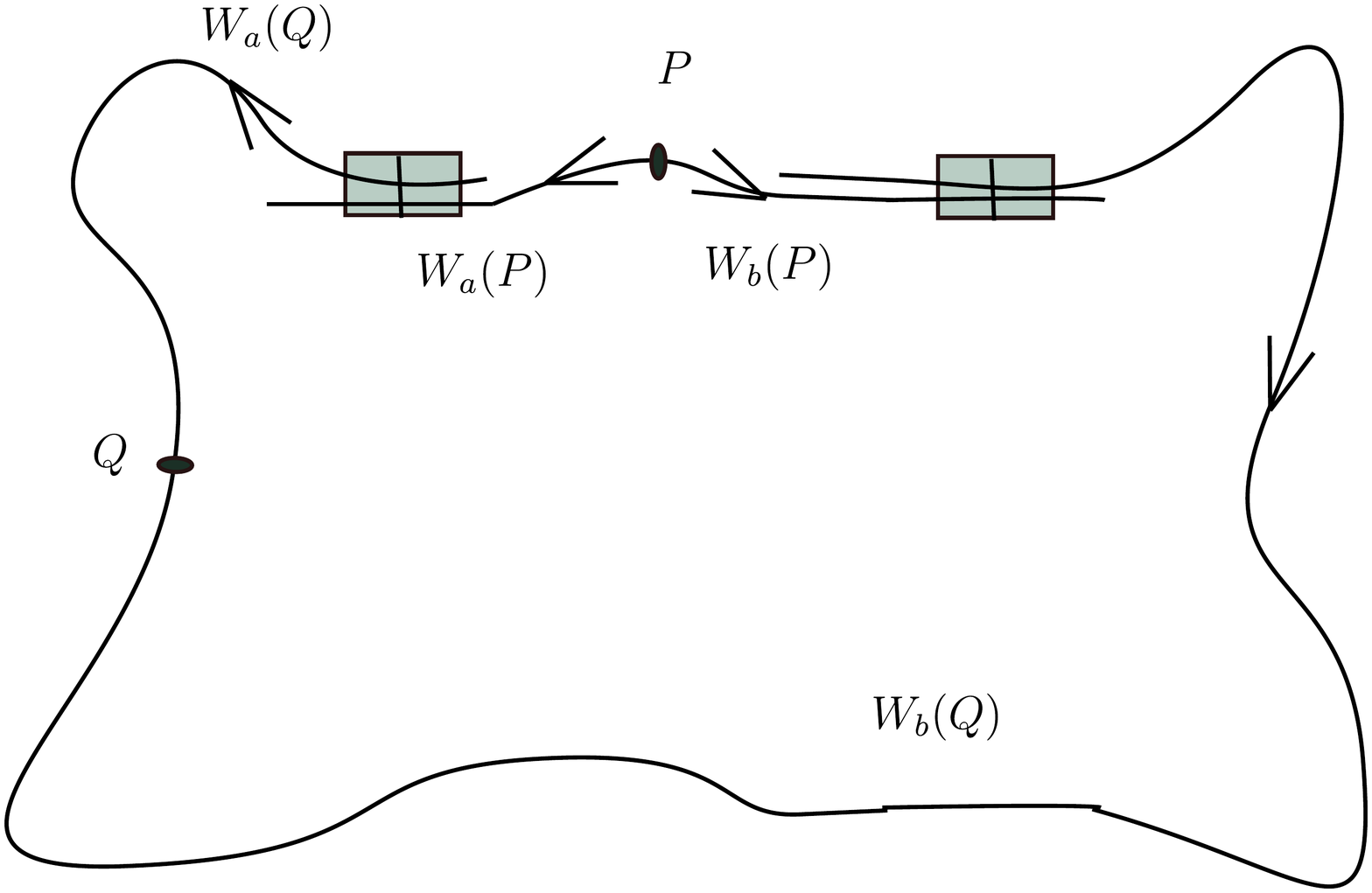}
 \caption{ Les coordonnées énergie-temps pour un bi-lien brisé.}\label{figbonnescoordonneeperturbees}
 \end{figure}
 
 \subsection{La restauration d'un bi-lien par des opérateurs} 
 Nous expliquons ici comment restaurer un bi-lien dans une bonne carte si on
 ne se soucie pas de l'impact sur les autres bi-liens.  Pour cela, on utilise
 la   carte énergie-temps construite dans le paragraphe précédent 
 ainsi que les notations introduites dans les paragraphes~\ref{secenergietemps} et~\ref{secbonnecarte}. Cette bonne carte est de classe $C^{r-1}$ et dépend de façon $C^1$ de la dynamique perturbée $F:\D\rightarrow \D$, elle-même proche de la dynamique non perturbée~$F_0$ qui possède un îlot stochastique. 
 
 Une bande verticale coupe le bi-lien de $F_0$  suivant deux composantes
 connexes, et donc une translation verticale $S_\psi: (x, y)\mapsto (x,
 y+\psi(x))$ à support dans une bande verticale~$V_a$ perturbe le bi-lien de
 $F_0$ ou  les arcs de branches qui lui correspondent pour $F$ en deux
 endroits\footnote{En réalité le bi-lien n'existe plus dans le cas perturbé,
   chaque branche du bi-lien de $F_0$ étant dédoublée par perturbation comme
   indiqué dans la figure  \ref{figbonnescoordonneeperturbees}. Aussi, couper $L_a$ ou $L_b$ dans le cas perturbé signifiera couper les deux branches correspondantes.}. Afin de pouvoir restaurer la première connexion~$L_a$ sans problème, Berger et Turaev   construisent  donc leur carte de manière à ce que $V_a$ ne  coupe $L_a$ qu'en un seul endroit.
 
 Comme la dynamique se lit au voisinage de $L_a\cap V_a$ comme la translation $(x, y)\mapsto (x-\tau, y)$, la partie $W^u(P)$ nous intéressant est le graphe d'une fonction $\tau$-périodique $w^u_a:[x_a-2\tau, x_a]\rightarrow \R$ et la partie de $W^s(Q)$ est aussi le graphe d'une  fonction $\tau$-périodique $w^s_a:[x_a-2\tau, x_a]\rightarrow \R$. L'écart entre les deux branches est alors la fonction $M_a(F)=\left( w^u_a-w^s_a\right)_{|[x_a-\tau, x_a]}$, fonction que nous voulons rendre nulle.  

 \begin{nota}
   Un nombre $\delta>0$ petit étant fixé, on note
   $C_\delta^{r-2}([x_a-2\tau, x_a], \R)$ l'ensemble des fonctions de classe
   $C^{r-2}$ définies sur $\R$ à support dans
   $[x_a-2\tau+\delta, x_a-\delta]$.
 \end{nota}

 Berger et Turaev montrent alors que l'opérateur $$ \psi\in C_\delta^{r-2}([x_a-2\tau, x_a], \R)\mapsto  M_a(S_\psi\circ F)$$ est de classe $C^1$, dépend continûment de $F$, et que dans  le cas non perturbé,   l'écart entre les deux branches est $M_a(S_\psi\circ F_0)(t)= M_a(F_0)+\psi(t)+\psi(t-\tau)$.

 Si on se restreint à  l'ensemble des fonctions $\psi=\rho_a\cdot\tilde\psi$ où $\rho_a\in C_\delta^\infty([x_a-2\tau, x_a], [0,1])$ vérifie $\rho_a(t)+\rho_a(t+\tau)=1$, alors l'opérateur $\Delta_a(F): \tilde\psi\mapsto M_a(S_{\rho_a\cdot\tilde\psi}\circ F)$ défini sur l'ensemble $\Pc^{r-2}_\tau$ des fonctions de classe $C^{r-2}$ et $\tau$-périodique de $\R$ vers~$\R$ est de classe~$C^1$ et est l'opérateur identité pour $F_0$\footnote{À condition d'identifier $\tilde\psi$ qui est $\tau$-périodique et $\tilde\psi_{|[x_a-\tau, x_a]}$.}. Aussi, ${\rm Id}_{\Pc^{r-2}_\tau} -\Delta_a(F_0)$ est  une contraction au voisinage de $0$. Pour $F$ assez proche de $F_0$, ${\rm Id}_{\Pc^{r-2}_\tau} -\Delta_a(F )$ est aussi une contraction sur ce voisinage, a donc un point fixe proche de $0$ et  $\Delta_a(F)$ a un zéro proche de $0$.  Ce zéro, que nous notons $\psi^a$,  rétablit la connexion~$L_a$.
 
 Il reste alors à rétablir la connexion $L_b$ à l'aide d'une translation verticale à support dans $V_b$, les branches correspondantes ayant été encore modifiées par la perturbation introduites pour $L_a$. La construction est plus compliquée car la bande $V_b$ coupe en deux composantes connexes le lien $L_b$ non perturbé. L'opérateur $M_b:F\mapsto \left( w^u_b-w^s_b\right)_{|[x_b, x_b+\tau ]}$ est construit comme l'était $M_a$ et $ \psi\in C_\delta^{r-4}([x_b, x_b+2\tau], \R)\mapsto  M_b(S_\psi\circ F)$ est de classe $C^1$, dépend continûment de $F$, et  dans  le cas non perturbé,   l'écart entre les deux branches devient alors 

 $M_b(S_\psi\circ F_0)(t)= M_b(F_0)(t)-\psi(t)-\psi(t+\tau)$
 $$+\frac{1}{2}\left( \psi(\frac{3x_b+\tau-x}{2})+\psi(\frac{3x_b+2\tau-x}{2})+\psi(\frac{3x_b+3\tau-x}{2})+\psi(\frac{3x_b+4\tau-x}{2})\right).$$
Une remarque fondamentale est alors qu'une fois que le lien $L_a$  a été restauré, $M_b(F)$ est d'intégrale nulle simplement parce que $F$ préserve l'aire. Aussi, on peut se restreindre aux fonctions d'intégrale nulle. Finalement, Berger et Turaev définissent sur l'espace~$\Pc^{r-4}_{\tau, 0}$ des fonctions $\psi:\R\rightarrow \R$   $\tau$-périodiques, de classe $C^{r-4}$  et d'intégrale nulle sur une période l'opérateur $\Delta_b(F): \tilde\psi\mapsto M_b(S_{\rho_b.\tilde\psi}\circ F)$ où cette fois le support de $\rho_b$ est dans $[x_b, x_b+2\tau]$ et on a toujours  $\rho_b(t)+\rho_b(t+\tau)=1$. Si on munit l'espace $\Pc^{r-4}_{\tau, 0}$ de la norme  
$$\| \psi\|=\max_{1\leq i\leq r-4}\| D^i\psi\|$$
alors ${\rm Id}_{\Pc^{r-4}_{\tau, 0}} -\Delta_b(F_0)$ devient encore une contraction et on conclut comme dans le cas précédent qu'il existe  $\psi^b$ qui  rétablit le lien.

Notons que comme cette nouvelle perturbation est à support dans $V_b$ qui ne rencontre pas $L_a$, cette nouvelle perturbation ne casse pas la partie restaurée $L_a$ et donc qu'on a ainsi restauré le bi-lien en entier.

 La perturbation globale  utilisée pour restaurer le bi-lien est finalement $S_\psi$ où $\psi=\rho_a.\psi^a+\rho_b.\psi^b$.

\subsection{De bonnes cartes}\label{secoplien}
\begin{figure}[H]
\centering \includegraphics[width=7.5cm]{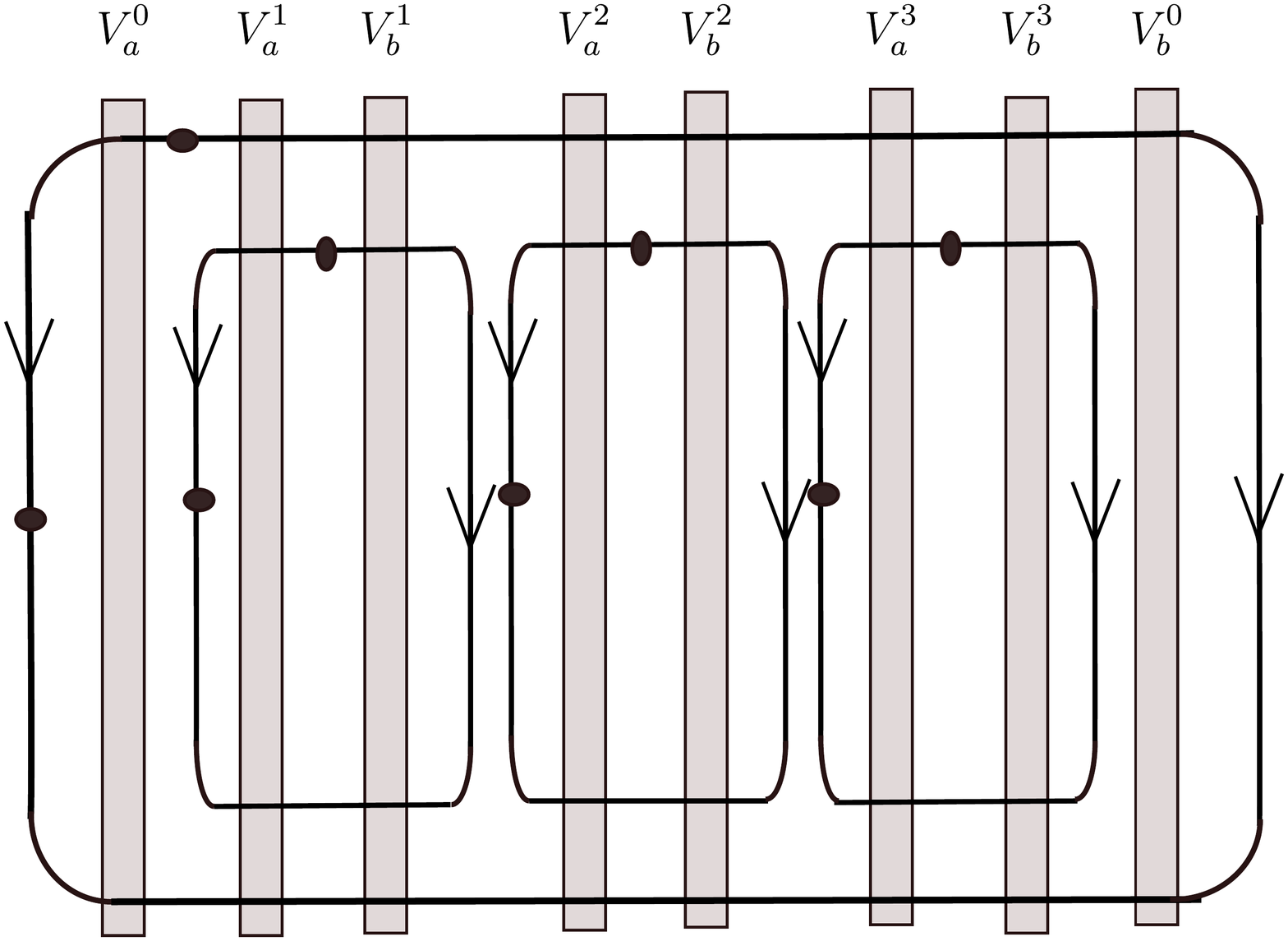}
 \caption{Une bonne carte pour l'îlot $\cI_0$.}\label{bonnecarte}
 \end{figure}

Il y a quatre bi-liens à restaurer dans le cas envisagé par Berger et
Turaev. Ils commencent par utiliser un difféomorphisme $\phi\in {\rm
  Diff}^r_\omega(\D)$ tel que $\hat f_0=\phi\circ f_0\circ \phi^{-1}$
détermine une bonne carte  pour chacun des quatre bi-liens avant perturbation
au sens 
du paragraphe~\ref{secbonnecarte}. La figure \ref{bonnecarte} résume ceci.

 Une carte analogue est construite dans le cas perturbé qui est une   carte  énergie-temps au sens 
 du paragraphe~\ref{secenergietemps} pour chaque bi-lien de l'îlot  perturbé.

 Il y a alors un ordre pour restaurer les bi-liens suivant la méthode que nous
 avons expliqué dans le paragraphe 
 précédent : tout d'abord, Berger et Turaev vont restaurer les bi-liens 1 à 3 en composant $\hat f=\phi\circ f\circ \phi^{-1}$ par une translation verticale de classe $C^{r-4}$ $(x, y)\mapsto (x, y+\psi_i(x))$   qui est à support dans $V_a^i\cup V_b^i$. Vu comment sont disposées les bandes verticales supports des perturbations, on peut restaurer les bi-liens 1  à  3 indépendamment,  mais cela va bien sûr changer les branches du dernier bi-lien 0.  Puis, on restaure ce lien 0 par une translation verticale $(x, y)\mapsto (x, y+\psi_0(x))$ de classe $C^{r-8}$ à support dans $V_a^0\cup V_b^0$, et ceci n'impacte pas les bi-liens 1 à  3 déjà restaurés car le support de la dernière perturbation ne les rencontre pas.
 
 La translation verticale $(x, y)\mapsto (x, y+\psi(x))$ de classe $C^{r-8}$ que l'on cherchait est donc définie par la  fonction $\displaystyle{ \psi=\sum_{i=0}^3\psi_ i}$.

\subsection{Synthèse} 
On commence par fixer un  voisinage de l'identité dans ${\rm
  Diff}_\omega^\infty(\D)$, qui correspond à  contrôler la distance $C^r$. On
construit alors $\hat f\in {\rm Diff}_\omega^\infty(\D)$ qui est $C^{r+8}$
proche de l'identité et qui admet comme application renormalisée une
application proche de l'îlot stochastique construit par Berger et Turaev dans
le paragraphe~\ref{secexilot}. Le corollaire~\ref{cororenorm} nous permet alors de rétablir les liens de l'îlot stochastique en utilisant les translations verticales de classe $C^r$ décrites 
dans le paragraphe~\ref{secoplien}. La proposition~\ref{Probustelien} nous dit alors que la dynamique obtenue, qui est de classe $C^r$, a un îlot stochastique et la proposition~\ref{proplissage} permet de rendre cet exemple de classe $C^\infty$.

\printbibliography


\end{document}
